\magnification=\magstep1
\baselineskip=13.9pt
\overfullrule 0pt

\centerline {\bf GROUP-LIKE ELEMENTS IN QUANTUM GROUPS}
\centerline {\bf AND FEIGIN'S CONJECTURE} 

\vskip .3cm 
\centerline {\bf Arkady Berenstein}
\centerline  {\it Department of Mathematics, Northeastern University, 
Boston}

\vskip .5cm 

\centerline  {\bf Abstract}
\item {}  Let $A$ be an arbitrary symmetrizable Cartan matrix of rank $r$,
and ${\bf n}={\bf n_+}$ 
be the standard maximal nilpotent subalgebra in the Kac-Moody
algebra associated 
with $A$ (thus, ${\bf n}$ is generated by $E_1,\ldots,E_r$ subject 
to the Serre relations).
Let $\hat U_q({\bf n})$ be the completion (with respect to 
the natural grading) of the quantized enveloping algebra 
of ${\bf n}$. For a sequence ${\bf i}=(i_1,\ldots,i_m)$ with $1\le i_k\le
r$, let $P_{\bf i}$ be a skew polynomial algebra generated by
$t_1,\ldots,t_m$ subject to the relations $t_lt_k=q^{C_{i_k,i_l}}t_kt_l$ 
($1\le
k<l\le m$) where $C=(C_{ij})=(d_ia_{ij})$ is the symmetric matrix
corresponding to $A$. We construct a group-like element ${\bf e}_{\bf i}\in
P_{\bf i}\bigotimes \hat U_q({\bf n})$. This element gives rise to the
evaluation homomorphism 
$\psi_{\bf i}:{\bf C}_q[N]\to P_{\bf i}$ given by $\psi_{\bf i}(x)=x({\bf
e}_{\bf i})$, where ${\bf C}_q[N]=U_q({\bf n})^0$ is the restricted dual of
$U_q({\bf n})$. Under a well-known isomorphism of algebras ${\bf C}_q[N]$
and $U_q({\bf n})$, the map $\psi_{\bf i}$ identifies with Feigin's
homomorphism $\Phi({\bf i}): U_q({\bf n})\to P_{\bf i}$. We prove that the
image of $\psi_{\bf i}$ generates the skew-field of fractions 
${\cal F}(P_{\bf i})$ if and only if ${\bf i}$ is a reduced expression of 
some
element $w$ in the Weyl group $W$; furthermore, in the latter case, ${\rm
Ker}~\psi_{\bf i}$ depends only on $w$ (so we denote 
$I_w:={\rm Ker}~\psi_{\bf
i}$). This
result generalizes the results in [5], [6] to the case of Kac-Moody
algebras. We also construct an  element 
${\cal R}_w\in \big({\bf C}_q[N]/I_w\big)\bigotimes \hat U_q({\bf n})$ which
specializes to  ${\bf e}_{\bf i}$ under the embedding  
${\bf C}_q[N]/I_w\hookrightarrow P_{\bf i}$.  The elements ${\cal R}_w$ are 
closely related to the quazi-$R$-matrix studied by G. Lusztig in [8].  
If ${\bf i},{\bf i}'$ are reduced expressions of the same element $w\in W$,
we have a natural isomorphism 
${\bf R}_{\bf i}^{\bf i'}:{\cal F}(P_{\bf i}) \to {\cal F}(P_{\bf i'})$ 
such that 
$({\bf R}_{\bf i}^{\bf i'}\otimes {\rm id})({\bf e}_{\bf i})={\bf e}_{\bf
i'}$. This leads to
identities between quantum exponentials. The maps ${\bf R}_{\bf i}^{\bf i'}$ 
are $q$-deformations of Lusztig's transition maps [8]. The existence of the
maps  ${\bf R}_{\bf i}^{\bf i'}$ leads to a  
surprising combinatorial corollary about skew-symmetric matrices associated  
with reduced expressions (cf. [12]).
  
\vfill
\eject
 
\noindent {\bf 0. Introduction and main results}

\medskip 

It is well-known that a quantum group is not a group. 
One of the goals of this chapter  is to introduce group-like elements for
quantum deformations of certain nilpotent algebraic groups.
In this section, we sketch our main results; more details will be given in
subsequent sections.

Consider  a maximal unipotent 
subgroup $N$ in a complex
simple algebraic group $G$. The group-like elements will be 
obtained
as quantum deformation of certain morphisms
$\pi_{\bf i}:{\bf C}^m\to N$ defined as follows.

Let $E_1,E_2,\ldots,E_r$ be standard generators of ${\bf n}$, 
the Lie algebra of $N$. For any sequence
${\bf i}=(i_1,\ldots,i_m)$ of indices (possibly with
repetitions), one defines a map
${\bf C}^m\to N$ by the formula 
$$\pi_{\bf i}(t_1,\ldots,t_m)= {\rm exp}(t_1E_{i_1})
{\rm exp}(t_2E_{i_2})\cdots
{\rm exp}(t_mE_{i_m})\eqno (0.1)$$
where ${\rm exp}: {\bf n}\to N$ is the exponential map. Note that $\pi_{\bf
i}$ is a
regular (algebraic) map.

It is well-known that $\pi_{\bf i}$ is
a birational isomorphism ${\bf C}^m\cong N$ if ${\bf i}=(i_1,\ldots,i_m)$
is a reduced 
expression of $w_0$, the longest element in the Weyl group $W$ of 
$G$. Furthermore, if ${\bf i}$ is a reduced expression  
of $w\in W$,  
then the closure in $N$ of the image of $\pi_{\bf i}$ depends only on $w$.

To construct a $q$-deformation of $\pi_{\bf i}$ we 
interpret the evaluation homomorphism $\pi_{\bf i}^*:{\bf C}[N]\to
{\bf C}[t_1,\ldots,t_m]$ as follows. 
First, we think of the product in (0.1) as an element 
$$\tilde\pi_{\bf i}\in {\bf C}[t_1,\ldots,t_m]\bigotimes 
\hat U({\bf n}),$$
where $\hat U({\bf n})$ is the completion of the universal 
enveloping algebra of ${\bf n}$ with respect to the natural grading. 

Second, ${\bf C}[N]$ can be identified with $U({\bf n})^0$, the 
restricted dual Hopf algebra. This gives rise to a natural pairing 
${\bf C}[N]\times \hat U({\bf n})\to {\bf C}$. Extending scalars 
from ${\bf C}$
to $P={\bf C}[t_1,\ldots,t_m]$ we see that each 
$f\in {\bf C}[N]$ becomes a linear form on $P\bigotimes 
\hat U({\bf n})$ with values in $P$. 
Then we have  
$$\pi_{\bf i}^*(f)=f(\tilde \pi_{\bf i}) \ . \eqno (0.2)$$ 

We construct the deformation of $\tilde \pi_{\bf i}$  in a more general
situation  when ${\bf n}$ is the standard maximal
nilpotent 
Lie subalgebra in a Kac-Moody algebra ${\bf g}$. Let us 
briefly introduce necessary definitions and notation.

\medskip

 Let $A=(a_{ij})$ be a symmetrizable Cartan
matrix of size $r \times r$. Denote by $C=(C_{ij})=(d_ia_{ij})$ 
the corresponding 
symmetric
matrix. Let ${\cal U}$ be the associative algebra over ${\bf
C}(q)$ generated by $E_1,\ldots,E_r$ subject 
to the {\it quantum Serre relations} (this is the quantized universal
enveloping 
algebra $U_q({\bf n})$ of the nilpotent part of the Kac-Moody algebra 
corresponding to $A$). 
The algebra ${\cal U}$ is graded by
${\bf Z}_+^r$ via ${\rm deg}(E_i)=\alpha_i$, the standard basis vector in
${\bf Z}_+^r$. Denote by $\hat {\cal U}$ the completion with respect 
to the grading. Following [8], Chapter 2, we consider ${\cal U}$ 
with the
structure of 
a {\it braided bialgebra} with the {\it braided} coproduct
$\Delta:{\cal U} \to {\cal U} \bigotimes {\cal U}$. 
Namely, $\Delta(E_i)=E_i\otimes 1 + 1\otimes E_i$, and 
$\Delta$ is a homomorphism of 
${\bf Z}_+^r$-graded algebras, where
the algebra structure on the tensor square of ${\cal U}$ differs 
from the
standard one by a twist (see Section 2 below for more details).  
It follows
that $\hat {\cal U}$ is a {\it complete bialgebra} with the coproduct 
$\hat \Delta: \hat {\cal U}\to \hat  {\cal U} \hat {\bigotimes} 
\hat {\cal U}$. 
The {\it quantum group} ${\cal A}$ 
is the restricted dual algebra of ${\cal U}$  
(if $A$ is of finite type then ${\cal A}$ can be identified with the
$q$-deformed ring of polynomial functions ${\bf C}_q[N]$). 
The natural evaluation 
pairing ${\cal A} \times \hat {\cal U}\to {\bf C}(q)$ 
will be  denoted  by $(x,E)\mapsto x(E)$.
 
\medskip 

Let ${\bf i}=(i_1,\ldots,i_m)$ be a sequence 
of integers with $1\le i_k\le r$. Denote by $P_{\bf i}$ 
the ${\bf C}(q)$-algebra generated by $t_1,\ldots,t_m$
subject to the relations
$$t_lt_k=q^{C_{i_k,i_l}}t_kt_l,~~~1\le k<l\le m \ . \eqno (0.3)$$
Let  $\hat {\cal U}_{\bf i}=P_{\bf i}\bigotimes \hat {\cal U}$ be 
the space of all series of the form 
$$\sum\limits_{\gamma \in {\bf Z}_+^r} t_\gamma \otimes
E_\gamma$$ 
where $t_\gamma\in P_{\bf i}$ and $ E_\gamma \in {\cal U}$ is a 
homogeneous
element of degree $\gamma$.  There 
is a standard algebra structure on $\hat {\cal U}_{\bf i}$. 
We identify 
$P_{\bf i}\bigotimes 1$ with $P_{\bf i}$ and   
$1\bigotimes \hat {\cal U}$ 
with $\hat {\cal U}$ so that $tE=Et=t\otimes E$ in 
$\hat {\cal U}_{\bf
i}$.  
Note also that $\hat {\cal U}_{\bf i}$ is a 
$P_{\bf i}$-bimodule in the 
standard way. The  
coproduct $\hat \Delta$ on $\hat {\cal U}$ extends naturally 
to the $P_{\bf i}$-bilinear map  
$$\hat \Delta_{\bf i}:\hat {\cal U}_{\bf i}\to \hat 
{\cal U}_{\bf i} \hat {\bigotimes
\limits_{P_{\bf i}}} \hat {\cal U}_{\bf i}$$ 
by the formula: $\hat \Delta_{\bf i}(\sum_\gamma t_\gamma 
E_\gamma)=\sum_\gamma
t_\gamma\Delta(E_\gamma)$. 
We call ${\bf e}\in \hat {\cal U}_{\bf i}$ 
a {\it group-like element} if $\hat \Delta_{\bf i}({\bf e})=
{\bf e}\otimes {\bf e}$. 
  
Finally, we define the
$q$-exponential 
$${\rm exp}_q(u)=\sum\limits_{n\ge 0} { u^n\over 
[n]_q!}\eqno (0.4)$$
where $[n]_q!=[1]_q[2]_q\cdots [n]_q$, $[l]_q=1+q+\cdots q^{(l-1)}$.

Now we can state our first main result.

\proclaim Theorem 0.1. For any sequence 
${\bf i}=(i_1,\ldots,i_m)$ and any
$c_1, \ldots, c_m \in {\bf C}(q)$, the
product 
$${\rm exp}_{q_{i_1}}(c_1 t_1E_{i_1})
{\rm exp}_{q_{i_2}}(c_2 t_2E_{i_2})\cdots
{\rm exp}_{q_{i_m}}(c_m t_mE_{i_m}) $$ 
is a {\rm group-like} element in $\hat {\cal U}_{\bf i}$, where
$q_i=q^{C_{ii}}$ for $i=1,\ldots,r$.

We prove Theorem 0.1 in Section 1 for more general {\it 
braided bialgebras}. 

We denote 
$${\bf e}_{\bf
i}={\rm exp}_{q_{i_1}}(t_1E_{i_1}){\rm exp}_{q_{i_2}}(t_2E_{i_2})\cdots
{\rm exp}_{q_{i_m}}(t_mE_{i_m}) \ .  \eqno (0.5)$$

As in the commutative case, we extend the evaluation pairing 
${\cal A}\times
\hat {\cal U} \to {\bf C}(q)$ to the $P_{\bf i}$-linear
pairing ${\cal A}\times \hat {\cal U}_{\bf i}\to P_{\bf i}$. As an
analogue of (0.2) we define the map $\psi_{\bf i}:
{\cal A}\to P_{\bf i}$
by 
$$\psi_{\bf i}(x):=x({\bf e}_{\bf i}) \ . \eqno (0.6)$$

\proclaim Corollary 0.2. The map $\psi_{\bf i}$ is an algebra homomorphism.

\noindent {\sl Proof}.  The definition of the pairing $(x,E)\to x(E)$ 
implies that  
$(xy)(u)=(x\otimes y)(\hat \Delta(u))$ for all $x,y\in {\cal A}$ and  
$u\in \hat {\cal U}_{\bf i}$, where $(x\otimes y)(u_1\otimes
u_2):=x(u_1)y(u_2)$ 
for any $u_1,u_2\in \hat {\cal U}_{\bf i}$. Thus, we have 
$$\psi_{\bf i}(xy)=(xy)({\bf e}_{\bf i})=
(x\otimes y)(\hat \Delta_{\bf i} ({\bf e}_{\bf i}))=
(x\otimes y)({\bf e}_{\bf
i}\otimes {\bf e}_{\bf i})=
x({\bf e}_{\bf i})y({\bf e}_{\bf i})=\psi_{\bf
i}(x)\psi_{\bf i}(y) \ .$$
Corollary 0.2 is proved. $\triangleleft$

\vskip .2cm

Expanding (0.5), we obtain the
following formula for $\psi_{\bf i}$: 
$$\psi_{\bf i}(x)=\sum\limits_{a_1,a_2\ldots,a_m\ge 0}
x\big(E_{i_1}^{[a_1]}E_{i_2}^{[a_2]}\cdots
E_{i_m}^{[a_m]}\big)t_1^{a_1}t_2^{a_2}\cdots t_m^{a_m}
\eqno (0.7)$$
where $E_i^{[n]}=\displaystyle{{1\over [n]_{q_i}!}E_i^n}$. Note 
that the sum in (0.7) is always finite. 

Under a well-known isomorphism  
${\cal A} \cong {\cal U}$ the homomorphism $\psi_{\bf i}$ becomes 
{\it Feigin's homomorphism}  $\Phi({\bf
i}):{\cal U}\to P_{\bf i}$ ([5] and Section~2 below).

\medskip 

Using the homomorphism $\psi_{\bf i}$, we can express the 
group-like element 
${\bf e}_{\bf i}$ in terms of the
{\it universal element} ${\cal R}\in {\cal A}
\bigotimes \hat {\cal U}$ defined as follows. Under the canonical 
isomorphism
between ${\cal A}
\bigotimes \hat {\cal U}$ and the space of linear maps ${\cal U}\to
\hat {\cal U}$, the element ${\cal R}$ corresponds to the 
inclusion ${\cal U}\hookrightarrow \hat {\cal U}$.
 
\proclaim Proposition 0.3. 
For any sequence ${\bf i}$ as above, we have 
$$(\psi_{\bf i}\otimes {\rm id})({\cal R})={\bf
e}_{\bf i} \ . \eqno (0.8)$$ 
The element ${\cal R}$ is uniquely determined by the equations  (0.8) 
for all ${\bf i}$.  

\noindent {\sl Proof}.  Let $B$ 
be a homogeneous basis in ${\cal U}$ (that is, $B$ is compatible 
with the ${\bf Z}_+^r$-grading),  
and $\{b^*\}$ be the dual basis in ${\cal
A}$ (so that $b^*(b')=\delta_{b,b'}$). By definition, 
$${\cal R}=\sum_{b\in B} b^*\otimes b \ . \eqno (0.9)$$ 
We have  
$$(\psi_{\bf i}\otimes {\rm id})({\cal R})=\sum_{b\in B} \psi_{\bf
i}(b^*)\otimes b=\sum_{b\in B} b^*({\bf e}_{\bf i})\otimes 
b={\bf e}_{\bf
i}$$
by definition (0.6) of $\psi_{\bf i}$  and by the formula 
$\sum_{b\in B}
b^*(u)\otimes b=u$ for any $u\in \hat {\cal U}_{\bf i}$.

It remains to check the uniqueness. Assume that there is 
another ${\cal R}'$
satisfying (0.8) for all ${\bf i}$. The equations (0.8) imply that  
$({\cal R}-{\cal R}')\in
{\rm Ker}~\psi_{\bf i}\bigotimes \hat {\cal U}$. By (0.7), an element 
$x\in {\cal A}$ is in the kernel of each $\psi_{\bf i}$ if and only if
$x$ vanishes at all monomials in $E_1,\ldots,E_r$. Hence,
${\cal R}'={\cal R}$, and we are done.  
$\triangleleft$

\vskip .1cm 

The element ${\cal R}$ was studied in  [8], Chapter 4  
in a slightly different setting; it is a version of the universal 
$R$-matrix for the ``braided" quantum double of ${\cal U}$. 

\vskip .2cm
 
Let $W$ be the Weyl group generated by simple
reflections $s_1,\ldots,s_r: {\bf Z}^r\to {\bf Z}^r$ defined by 
$s_i(\alpha_j)=\alpha_i-a_{ij}\alpha_i$ for all $i,j$. We say 
that ${\bf
i}=(i_1,\ldots,i_m)$ is
a reduced expression of $w\in W$ if $w=s_{i_1}\cdots s_{i_m}$ and 
this factorization of $w$ is the shortest possible. We denote by 
$R(w)$ the set of all reduced
expressions of $w$. We also reserve notation $w_0$ for the
longest element in $W$ if $W$ is finite.

For a sequence ${\bf i}=(i_1,\ldots,i_m)$ let ${\cal U}({\bf i})$ be
the subspace in ${\cal U}$ spanned by all monomials
$E_{i_1}^{a_1}E_{i_2}^{a_2}\cdots
E_{i_m}^{a_m}$. It is-well known ([7], Section 4.4, or [9]) that if 
${\bf i}\in
R(w)$ then ${\cal U}({\bf i})$ depends only on $w$. So we denote 
${\cal
U}(w):={\cal U}({\bf i})$ for all ${\bf i}\in R(w)$. It is also 
well-known
that ${\cal U}({w_0})={\cal U}$ if $W$ is finite. Further, if ${\bf i}$
is not reduced 
then there is a subsequence ${\bf i}'$ of ${\bf i}$ such that 
${\bf i}'$ is
a reduced expression and ${\cal U}({\bf i})={\cal U}({\bf i'})$.  
For the convenience of the reader we will prove these assertions in 
Section 2. 

Now we can give a complete description of  ${\rm Ker}~\psi_{\bf i}$ 
using
(0.7)
and the above discussion.

\proclaim Lemma 0.4. The kernel of $\psi_{\bf i}$ is the
orthogonal complement of ${\cal U}({\bf i})$:
$${\rm Ker}~ \psi_{\bf i}=\{x\in {\cal A}:
x(u)= 0 ~ {\rm for~ all} ~ u\in {\cal U}({\bf i})\} \ .$$
Furthermore, 
\item {(i)} If ${\bf i}$ is not reduced then there is a reduced 
subsequence
${\bf i}'$
of ${\bf i}$ such that ${\rm Ker}~ \psi_{\bf i}={\rm Ker}~ 
\psi_{\bf i'}$.
\item {(ii)} For every $w\in W$ and ${\bf i},{\bf i'}\in R(w)$ 
we have ${\rm Ker}~
\psi_{\bf i}={\rm Ker}~ \psi_{\bf i'}$.
\item {(iii)} If $W$ is finite then ${\rm Ker}~ \psi_{\bf i}=\{0\}$ for 
any ${\bf i}\in R(w_0)$. 

Denote $I_w:={\rm Ker}~ \psi_{\bf i}$ for ${\bf i}\in R(w)$. Our 
next main result is the following.   

\proclaim Theorem 0.5. 
For every $w\in W$ and 
${\bf i} \in R(w)$, the image of $\psi_{\bf i}$ generates 
the (skew) field of fractions of $P_{\bf i}$. Hence, $\psi_{\bf i}$ 
induces an isomorphism 
of the fields of fractions 
$$\overline \psi_{\bf i}:{\cal F}\big({\cal A}/I_w\big)\widetilde \to
{\cal F}\big(P_{\bf i}\big) \ . \eqno (0.10)$$
In particular, if $W$ is finite and $w = w_0$ then $\overline \psi_{\bf i}$
is an 
isomorphism between ${\cal F}\big({\cal A}\big)$ and ${\cal F}\big(P_{\bf
i}\big)$. 

Here the symbol ${\cal F}$ stands for the skew-field of fractions. 
 We review the necessary definitions and results in the appendix below. 

\medskip 
  
The last statement in Theorem 0.5 coincides with Feigin's conjecture 
([5],~[6]) 
(stated for ${\cal U}$ rather than for ${\cal A}$).  The
conjecture was proved in [5] for the type $A_r$ and $w=w_0$ by a
direct computation involving some specific reduced word ${\bf i}\in R(w_0)$. 
The conjecture was further
generalized by A. Joseph ([6]) to any $w\in W$ in the assumption 
that $W$ is finite, and proved  by  geometric arguments.

We give an algebraic proof of Theorem 0.5 in Section 2 without the  
assumption that $W$ is finite. The following proposition plays the 
crucial role in the proof. 
 
\proclaim Proposition 0.6. For every reduced expression ${\bf i}$ of 
some $w\in W$, there is an element
$x=x({\bf i})\in {\cal A}$ such that 
$\psi_{\bf i}(x)=t_1^{a_1}t_2^{a_2}\cdots t_m^{a_m}$
with $a_1>0$. 

We prove Proposition 0.6 in Section 3.

\proclaim Corollary 0.7. The skew-field ${\cal F}\big 
(\psi_{\bf i}({\cal A})\big)$ coincides with ${\cal F}(P_{\bf i})$ 
if and only if ${\bf i}$ is a reduced expression. 

\noindent {\sl Proof}. We need to prove the ``only if" part
(the ``if" part is the assertion of Theorem 0.5). Assume that 
${\bf i}=(i_1,\ldots,i_m)$ is not reduced. Then, by Lemma 0.4(i) and 
Theorem 0.5, there is a reduced subsequence ${\bf i'}=
(i_{k_1},\ldots,i_{k_l})$ of ${\bf i}$ such that ${\rm
Im}~\psi_{\bf i'}\cong {\rm Im}~\psi_{\bf i}$. It follows that
${\cal F}({\rm Im}~\psi_{\bf i})\cong {\cal
F}({\rm Im}~\psi_{\bf i'})\cong {\cal F}(P_{\bf i'})$. 
But ${\cal F}(P_{\bf i'})\not \cong {\cal F}(P_{\bf i})$ 
since these skew-fields have different Gel'fand-Kirillov dimensions 
(see [12], Proposition  2.18). Corollary 0.7 is proved. $\triangleleft$ 

\medskip

Let us illustrate the above results in the case
when $A$ is of type $A_{n-1}$.
In this case, one can show that ${\cal A}$ is generated 
by the elements
$x_{ij}$ with $1 \le i < j \le n$ subject to 
the following relations (cf. [2],~[4]):
$$x_{ij}x_{kl}=x_{kl}x_{ij},~x_{il}x_{jk}=x_{jk}x_{il} \quad
(1\le i<j<k<l\le n) \ ,\eqno (0.11)$$
$$\displaystyle{x_{ik}={qx_{ij}x_{jk}-x_{jk}x_{ij}\over
q-q^{-1}}},~~x_{ij}x_{ik}=qx_{ik}x_{ij},~~x_{jk}x_{ik}
=q^{-1}x_{ik}x_{jk} \,\, 
(1\le i<j<k\le n) \ . $$
The elements $x_{ij}$ are $q$-deformations of the matrix entries 
considered as polynomial functions on the group $N$ of 
the unipotent upper-triangular
matrices. 
So we arrange the $x_{ij}$ into the matrix 
$X = I+\sum\limits_{i<j}
x_{ij} E_{ij}$, where $I$ is the identity matrix, and 
the $E_{ij}$ are the matrix units.
Let $\psi_{\bf i}(X)$ be the
$n\times n$-matrix (over $P_{\bf i}$)
obtained from $X$ by applying  
$\psi_{\bf i}$ to each matrix entry. 
The following proposition will be proved in Section~2. 

\proclaim Proposition 0.8. For any sequence ${\bf i} = 
(i_1, \ldots, i_m)$,
the matrix $\psi_{\bf i}(X)$ admits
the following matrix factorization 
$$\psi_{\bf i}(X)=(I+t_1E_{i_1,i_1+1})(I+t_2E_{i_2,i_2+1})\cdots
(I+t_mE_{i_m,i_m+1}) \ . \eqno (0.12)$$

If ${\bf i} \in R(w_0)$	then, applying the inverse isomorphism 
$(\overline \psi_{\bf i})^{-1}: {\cal F}\big(P_{\bf i}\big) \to
{\cal F}\big({\cal A}\big)$ to the factorization (0.12), we 
obtain the factorization of the matrix $X$ over 
${\cal F}\big({\cal A}\big)$:
$$X=(I+\tilde t_1E_{i_1,i_1+1})(I+\tilde t_2E_{i_2,i_2+1})\cdots
(I+\tilde t_mE_{i_m,i_m+1}) \ ,$$
where $\tilde t_k = (\overline \psi_{\bf i})^{-1} (t_k)$.
This factorization is a $q$-deformation of the one studied in [1];
it can be shown that such a factorization is unique. 
The explicit formulas for $\tilde t_k$ in terms of the matrix entries 
$x_{ij}$
will be given in a separate publication. The above
factorizations of the matrix $X$ (and, more generally, ${\cal R}$ for
quantum groups of finite type) were studied in [11]. 

 \medskip

Let us return to the general situation and discuss some corollaries of 
Theorem 0.5. 
For every $w\in W$ and ${\bf i},{\bf i'}\in R(w)$ there is an 
isomorphism
of skew fields 
$${\bf R}_{\bf i}^{\bf i'}:{\cal F}\big(P_{\bf i})\widetilde \to {\cal
F}\big(P_{\bf i'}\big)$$  
defined by ${\bf R}_{\bf i}^{\bf i'}
=\overline  \psi_{\bf i'}\circ (\overline \psi_{\bf i})^{-1}$.  
We extend
it to the isomorphism ${\bf R}_{\bf i}^{\bf i'}\otimes {\rm id}: 
{\cal
F}\big(P_{\bf i}\big)\bigotimes \hat {\cal U}\to {\cal
F}\big(P_{\bf i'}\big)\bigotimes \hat {\cal U}$.  
In the following proposition, every element ${\bf e}_{\bf i}$ given by 
(0.5) is regarded as 
an element of ${\cal F}\big(P_{\bf i}\big)\bigotimes \hat {\cal U}$.  

\proclaim Proposition 0.9. For every ${\bf i},{\bf i'}\in R(w)$, we
have 
$({\bf R}_{\bf i}^{\bf i'}\otimes 
{\rm id})({\bf e}_{\bf i})={\bf e}_{\bf i'}$.

\noindent {\sl Proof}. Let ${\bf p}_w:{\cal A}\to {\cal A}/I_w$ 
be the canonical 
projection. Denote ${\cal R}_w=({\bf p}_w\otimes {\rm id})({\cal R})$. 
Note that, similarly to ${\cal R}$, the element 
${\cal R}_w\in {\cal A}/I_w\bigotimes \hat {\cal U}$  
corresponds to the inclusion 
${\cal U}(w)\hookrightarrow {\cal U}$.  
Proposition 0.3 implies that 
$$(\overline \psi_{\bf i} \otimes {\rm id})({\cal R}_w)={\bf e}_{\bf
i}\eqno (0.13)$$
for every ${\bf i}\in R(w)$. We are done since 
${\bf R}_{\bf i}^{\bf i'}\otimes {\rm id}=
(\overline  \psi_{\bf i'}\otimes {\rm id})
\circ (\overline \psi_{\bf i}\otimes {\rm id})^{-1}$. 
$\triangleleft$
 
\medskip

Proposition 0.9 implies some  identities between 
quantum exponentials. 
For two reduced expressions ${\bf i}=(i_1,\ldots,i_m)$ and
${\bf i}'=(i'_1,\ldots,i'_m)$ of an element 
$w\in W$, let $t_1,\ldots,t_m$ 
(resp. $t'_1,\ldots,t'_m$) be the standard generators 
of $P_{\bf i}$ (resp. $P_{\bf i'}$).

\proclaim Corollary  0.10. 
The following identity holds in the algebra 
${\cal F}(P_{\bf i})\bigotimes \hat {\cal U}$:
$${\rm exp}_{q_{i_1}}(t_1E_{i_1})\cdots 
{\rm exp}_{q_{i_m}}(t_mE_{i_m})
={\rm exp}_{q_{i'_1}}(p_1E_{i'_1})\cdots {\rm exp}_{q_{i'_m}}
(p_mE_{i'_m}) \ , \eqno (0.14)$$
where $p_k={\bf R}_{\bf i'}^{\bf i}(t'_k)$ for $k=1,\ldots,m$. 
Identity (0.14) remains true under
the rescaling $E_i\mapsto c_iE_i$ for any  
$c_i\in {\bf C}(q)$, $i=1,\ldots,r$.

\noindent {\bf Example 1}. Let $A$ be the Cartan 
matrix of type $A_2$ . Then the Weyl group $W$ is the symmetric 
group $S_3$, and 
$w_0=s_1s_2s_1=s_2s_1s_2$. Take  ${\bf
i}=(121), \, {\bf i}'=(212)$, and denote by $t_1,t_2,t_3$ and 
$t'_1,t'_2,t'_3$ the generators of 
$P_{(121)}$ and $P_{(212)}$ respectively.  Then 
${\bf R}_{(212)}^{(121)}(t_k')=p_k$ for $k=1,2,3$, where 
$$p_1=t_2t_3(t_1+t_3)^{-1},~p_2=t_1+t_3,~p_3=(t_1+t_3)^{-1}t_1t_2 \ . $$
This is a consequence of the matrix equation 
$$(I+p_1E_{23})(I+p_2E_{12})(I+p_3E_{23})
=(I+t_1E_{12})(I+t_2E_{23})(I+t_3E_{12}) \ ,$$
which follows from the factorization (0.12). 

The identity (0.14)  takes the form
$${\rm exp}_{q^2}(c_1t_1E_1)
{\rm exp}_{q^2}(c_2t_2E_2){\rm exp}_{q^2}(c_1t_3E_1)=
{\rm exp}_{q^2}(c_2p_1E_2)
{\rm exp}_{q^2}(c_1p_2E_1){\rm exp}_{q^2}(c_2p_3E_2) \eqno (0.15)$$
for any $c_1,c_2\in {\bf C}(q)$.
Expanding both sides of (0.15) and comparing the components of degree  
$2 \alpha_1 + \alpha_2$, we obtain the quantum Serre relation 
$$E_1^2E_2-(q+q^{-1})E_1
E_2E_1-E_2E_1^2= 0 \ .$$ 
We also note that setting $c_1=1$ and $c_2=0$ in (0.15) yields the 
familiar rule
$${\rm exp}_{q^2}(t_1E_1){\rm exp}_{q^2}(t_3E_1)
={\rm exp}_{q^2}((t_1+t_3)E_1) \ . $$

The identity (0.15) appeared in [11], Section 10.4; it was proved there  
by a straightforward computation. 

\medskip 

We conclude the introduction by a surprising combinatorial 
consequence of the above results.
To a sequence ${\bf i}=(i_1,\ldots,i_m)$ we associate a 
skew-symmetric $m\times m$-matrix $S({\bf i})$ by the formula:
$$S({\bf i})=\sum_{1\le k<l \le m} C_{i_k,i_l}(E_{kl}-E_{lk}) \ . 
\eqno (0.16)$$ 
We say that two $m\times m$ matrices $S$ and $S'$ are 
{\it equivalent} if there 
is a matrix $T\in SL_m({\bf Z})$ such that $S'=TS T^t$  
(where $T^t$ is the transpose of $T$). 

\proclaim Proposition 0.11. For every $w$ and ${\bf i},
{\bf i'}\in R(w)$,
the matrices $S({\bf i})$ and
$S({\bf i'})$ are equivalent.

This follows from the fact that the skew-fields 
${\cal F}(P_{\bf i})$ and
${\cal F}(P_{\bf i'})$ are isomorphic, in view of a 
general result by 
A.~Panov [12] (see also Section 2 below). 
Proposition~0.11 essentially says that 
there exists an isomorphism 
${\cal F}(P_{\bf i'})\to {\cal F}(P_{\bf
i})$ which takes every generator $t_k'$ to a monomial in 
$t_1,\ldots,t_m$.
(Note that the isomorphism ${\bf R}_{\bf i'}^{\bf i}$ 
considered above, 
in general does not have this property). 

\medskip

The material is organized as follows.
In Section~1 we introduce braided bialgebras and prove some results
about them, including the generalization of Theorem~0.1. 
The quantum group ${\cal A}$ associated with a symmetrizable 
Cartan matrix
is studied in Section~2, which contains the proofs of Lemma~0.4, 
Theorem~0.5 
(modulo Proposition~0.6), and Propositions~0.8 and 0.11.
Section~3 is devoted to the proof of Proposition~0.6; our proof 
is based 
on the properties
of extremal vectors in simple $U_q({\bf g})$-modules. 
In Appendix we review necessary definitions and results about 
non-commutative fields of 
fractions. 
 
\medskip

\noindent {\bf Acknowledgments}. The work presents results first 
appeared in my Ph.D.  dissertation (Northeastern University, May 1996). 
I am very grateful to my advisor  
Professor
Andrei Zelevinsky, for the continual support and help 
throughout the writing of this paper. This work benefited greatly 
from my conversations with J. Bernstein,
G. Lusztig, S. Majid, A. Postnikov, Y. Soibelman, and J. Towber. 
I am especially indebted to Professor David Kazhdan of Harvard, first for
making me welcome at Harvard, and second  
for teaching me geometric and categorical approach to the subject. 
Thanks are due the faculty of Department of Mathematics at 
Northeastern University, for the help during my stay there.

\bigskip

\noindent  {\bf 1. Results on braided bialgebras}

\medskip

Let ${\bf k}$ be a field and ${\cal U}$ be a 
${\bf Z}_+^r$-graded ${\bf k}$-algebra: ${\cal U}
=\bigoplus\limits{\cal U}(\gamma)$, 
the sum over $\gamma\in {\bf Z}_+^r$. 
We assume that ${\cal U}(0)={\bf k}$ and every ${\cal U}(\gamma)$ is 
finite-dimensional. 
Let ${\bf q}=(q_{ij})$, $1\le i,j\le r$ be a  
$r\times r$-matrix with all
$q_{ij}\in {\bf k},q_{ij}\ne 0$. Following G. 
Lusztig ([8]) we associate with 
${\bf q}$ an algebra structure on the vector space  
${\cal U}\bigotimes
{\cal U}$. For any two homogeneous elements 
$b\in {\cal U}(m_1,\ldots,m_r)$
and $c\in {\cal U}(n_1,\ldots,n_r)$ we set 
$$Q(b,c)=Q((m_1,\ldots,m_r),(n_1,\ldots,n_r))=
\prod_{i,j=1}^r q_{ij}^{m_in_j} \ . \eqno (1.1)$$  
We define the ${\bf q}$-{\it braided} 
multiplication in ${\cal U}\bigotimes {\cal U}$ by 
$$(a\otimes b)(c\otimes d):=Q(b,c)(ac\otimes bd)\eqno (1.2)$$
for any   homogeneous elements $b,c$ of ${\cal
U}$ and any $a,d\in {\cal U}$.

It is easy to see that (1.2) makes 
${\cal U}\bigotimes {\cal U}$ into a 
${\bf Z}_+^r$-graded associative algebra 
(with the standard grading $({\cal
U}\bigotimes_{\bf k} {\cal U})(\gamma)=\bigoplus_{\gamma'}~ {\cal
U}(\gamma')\bigotimes {\cal U}(\gamma-\gamma')$). 
This algebra will be denoted  by ${\cal U}
{\bigotimes}_{\bf q}{\cal U}$ and called the 
${\bf q}$-{\it braided} tensor square of ${\cal U}$. 

\vskip .2cm 

We call ${\cal U}$   
a ${\bf q}$-{\it braided bialgebra} if 
\item {(i)} there is a homomorphism of ${\bf Z}_+^r$-graded algebras
$\Delta:{\cal U}\to {\cal U} {\bigotimes}_{\bf q}
{\cal U}$ satisfying the coassociativity constrain 
(we call $\Delta$ the {\it coproduct});
\item {(ii)} There is a {\it counit} homomorphism of algebras 
$\varepsilon:{\cal U}\to
{\cal U}(0)={\bf k}$ satisfying
$$(\varepsilon\otimes {\rm id})\circ \Delta=({\rm id}\otimes
\varepsilon)\circ \Delta={\rm id}, ~\varepsilon(1)=1.\eqno (1.3)$$

This definition implies  that for every $u\in {\cal U}$, 
$$\Delta(u)=u\otimes 1+1\otimes u+\sum_n u_n\otimes u'_n\eqno( 1.4)$$
where all $u_n,u'_n$ are homogeneous elements of nonzero degrees. In
particular,
$\Delta(u)=u\otimes 1+1\otimes u$ for every $x\in {\cal
U}(\alpha_1)\bigoplus \cdots \bigoplus {\cal
U}(\alpha_r)$ where $\alpha_1,\ldots,\alpha_r$ is the standard 
basis in ${\bf Z}_+^r$. Another  consequence of this  
definition is that $\varepsilon(x)=0$
for any $x\in {\cal U}(\gamma),\gamma\ne 0$. 

Note that the algebra ${\cal U}$ from the introduction, 
associated to a symmetrizable
Cartan matrix $A$, is a ${\bf q}$-braided algebra, where 
$q_{ij} = q^{C_{ij}}$.
Another example is the {\it free} algebra generated by 
$E_1, \ldots, E_r$,
where ${\bf q}$ is arbitrary.

Let $\hat {\cal U}$ be
the completion of ${\cal U}$ with respect to the grading, 
that is, the space
of all formal series $\hat u=\sum\limits_{\gamma\in {\bf Z}_+^r} 
u_\gamma$,
where $u_\gamma\in {\cal U}(\gamma)$.  
Clearly, $\hat {\cal U}$ is an algebra. 
The coproduct in ${\cal U}$ extends to  
$\hat \Delta: \hat {\cal U}\to \hat {\cal U}\hat 
{\bigotimes}_{\bf q} 
\hat {\cal U}$ 
so $\hat {\cal U}$ becomes a {\it complete bialgebra}.

Now we fix a positive integer $m$ and consider a sequence ${\bf
i}=(i_1,i_2,\ldots,i_m)$ of integers with  
$1\le i_k\le r$. Let ${\bf q}=(q_{ij})$
be the matrix used in the definition of ${\cal U}$. Consider a ${\bf
k}$-algebra $P_{\bf i}=P_{\bf i,q}$  generated by 
$t_1,\ldots,t_m$ subject
to the following relations:
$$t_lt_k=q_{i_k,i_l}t_kt_l\eqno (1.5)$$ for all $1\le k<l\le m$. 
  
Define $\hat {\cal U}_{\bf i}=P_{\bf i}\bigotimes_{\bf k} 
\hat {\cal U}$, the space of formal series of the  
form $\sum_\gamma t_\gamma 
\otimes u_\gamma$, 
where $t_\gamma\in P_{\bf i}$ and $u_\gamma\in {\cal U}(\gamma)$.
We consider $\hat {\cal U}_{\bf i}$ with the standard 
algebra structure (so we can write $tu=ut=t\otimes u$).

Consider the completed tensor square 
${\cal V}_{\bf i}=\hat {\cal U}_{\bf
i}\hat {\bigotimes\limits_{P_{\bf i}}} \hat {\cal U}_{\bf i}$ 
where the left factor is regarded as a 
right $P_{\bf i}$-module and the right 
factor as a left
$P_{\bf i}$-module. Note that ${\cal V}_{\bf i}$ is a 
$P_{\bf i}$-bimodule. In ${\cal V}_{\bf i}$, we can write 
$t(u\otimes v)=(tu)\otimes v=u\otimes (tv)=(u\otimes v)t$
for any $u,v\in {\cal U}, t\in P_{\bf i}$. 
Under the standard identification 
${\cal V}_{\bf i}\cong P_{\bf i} \bigotimes  
\hat {\cal U}\hat {\bigotimes}_{\bf q} \hat {\cal U}$ 
this bimodule ${\cal V}_{\bf i}$ becomes an algebra.

There is a  natural morphism of $P_{\bf i}$-bimodules 
$$\hat \Delta_{\bf i}:\hat {\cal U}_{\bf i}\to {\cal V}_{\bf i}$$
which is the $P_{\bf i}$-linear extension of 
the coproduct $\hat \Delta$ on  $\hat {\cal U}$. Clearly, 
$\hat \Delta_{\bf i}$ is an algebra homomorphism.

Let ${\bf E}=(E_1,\ldots,E_m)$ be the 
family of elements $E_k\in {\cal
U}(\alpha_{i_k})$. We define 
an element ${\bf e}_{\bf i} = {\bf e}_{\bf i,E} 
\in \hat {\cal U}_{\bf i}$ as follows:
$${\bf e}_{\bf i,E}={\rm exp}_{q_1}(t_1E_1)
{\rm exp}_{q_2}(t_2E_2)\cdots
{\rm exp}_{q_m}(t_mE_m)\eqno (1.6)$$  
where $q_k=q_{i_k,i_k}$ for $k=1,\ldots,m$, 
and ${\rm exp}_{q_k}$ stands for 
the quantum exponential defined by (0.4). 

The following result extends Theorem~0.1 to arbitrary 
${\bf q}$-braided algebras. 

\proclaim Theorem 1.1. For any sequence ${\bf i}$ and any 
family ${\bf E}=(E_k)$ as above the element 
${\bf e}_{\bf i}={\bf e}_{\bf i,E}$ is 
a {\rm group-like} element in $\hat {\cal U}_{\bf i}$, 
i.e., $\hat \Delta_{\bf i}({\bf e}_{\bf i})=
{\bf e}_{\bf i}\otimes {\bf
e}_{\bf i}$.

\noindent {\sl Proof }. We need the following.  
\proclaim Lemma 1.2. \item {(a)} Each factor ${\bf e}_k=
{\rm exp}_{q_k}(t_kE_k)$ of ${\bf
e}_{\bf i}$ is a group-like element in $\hat {\cal U}_{\bf i}$. 
\item {(b)} $(1\otimes {\bf e}_k)({\bf e}_l
\otimes 1)=({\bf e}_l\otimes 1) (1 \otimes
{\bf e}_k)$ for any $1\le k<l\le m$.

\noindent {\sl Proof}. (a) Denote $E=t_kE_k$. Since 
$\Delta(E_k)=
E_k\otimes 1+1\otimes E_k$, for each $k$ we have 
$$\hat \Delta_{\bf i}(E)=t_k(E_k\otimes 1+1\otimes E_k)=
E\otimes 1+1\otimes E.$$ 
Denote $x=E\otimes 1, y=1\otimes E$. Let us show that that 
$yx=qxy$ where
$q:=q_k=q_{i_k,i_k}$. 
Indeed, 
$$yx=(1\otimes E)(E\otimes 1)=(1\otimes t_kE_k)(t_kE_k\otimes 1)=
t_k^2(1\otimes E_k)(E_k\otimes 1)$$
$$=Q(E_k,E_k)t_k^2(E_k\otimes E_k)=
q_{i_k,i_k}(t_kE_k\otimes t_kE_k)=qxy \ .$$ 
Further, we obtain  
$$\hat \Delta_{\bf i}({\rm exp}_q(E))
={\rm exp}_q(\hat \Delta_{\bf i}(E))
={\rm exp}_q(x+y)$$ and 
$${\rm exp}_q(E)\otimes
{\rm exp}_q(E)=({\rm exp}_q(E)\otimes 1)(1\otimes {\rm exp}_q(E))$$
$$=({\rm exp}_q(E\otimes 1))({\rm exp}_q(1\otimes
E))={\rm exp}_q(x) {\rm exp}_q(y) \ . $$ Then the well-known rule 
for the quantum exponentials. 
$${\rm exp}_q(x+y)={\rm exp}_q(x) {\rm exp}_q(y)$$
(provided that $yx=qxy$) implies that 
$\hat \Delta_{\bf i}({\rm exp}_q(E))=
{\rm exp}_q(E)\otimes {\rm exp}_q(E)$. Part (a) is proved. 

\smallskip

(b)  Denote $E=t_kE_k$ and $E'=t_lE_l$. By definition of ${\cal
U} {\bigotimes}_{\bf q} {\cal U}$,
$$(1\otimes E)(E'\otimes 1)=t_kt_l(1\otimes E_k)(E_l\otimes
1)=q_{i_k,i_l}t_kt_l(E_l\otimes E_k) \ . $$
The commutation relations (1.5) imply that 
$$(1\otimes E)(E'\otimes 1)=t_kt_l(1\otimes E_k)(E_l\otimes 1)
=t_lt_k(E_l\otimes E_k)=E'\otimes E=(E'\otimes 1) (1\otimes E).$$
It follows that $(1\otimes f(E))(g(E')\otimes 1)=f(E')\otimes 
g(E)=(f(E')\otimes 1) (1\otimes g(E))$
for any polynomials $f$ and $g$ in one variable. Passing to the 
completion, we see that 
$f$ and $g$ can also be power series in the above formula.  
Taking  $f(E):={\bf e}_k={\rm exp}_{q_k}(E)$ and $g(E'):={\bf e}_l={\rm
exp}_{q_l}(E')$ completes the proof of part (b).  
Lemma 1.2 is proved.   
$\triangleleft$

\medskip 

We are ready to complete the proof of Theorem 1.1 now. 
Recall that we use the
shorthand ${\bf e}_k={\rm exp}_{q_k}(t_kE_k)$ so 
${\bf e}_{\bf i}={\bf
e}_1{\bf e}_2\cdots {\bf e}_m$. 

Using Lemma 1.2 and the fact that 
$(a\otimes 1)(1\otimes b)=
a\otimes b$ for any  $a,b\in \hat {\cal U}_{\bf i}$,
we obtain 
$$\hat \Delta_{\bf i}({\bf e}_{\bf i})=
\hat \Delta_{\bf i}({\bf e}_1{\bf
e}_2\cdots {\bf e}_m)=
\hat \Delta_{\bf i}({\bf e}_1)
\hat \Delta_{\bf i}({\bf e}_2)\cdots \hat
\Delta_{\bf i}({\bf e}_m)=
({\bf e}_1\otimes {\bf e}_1)({\bf e}_2\otimes {\bf e}_2)\cdots 
({\bf e}_m\otimes {\bf e}_m)$$
$$=({\bf e}_1\otimes 1)(1\otimes {\bf e}_1)({\bf e}_2\otimes 1)
(1\otimes
{\bf e}_2)\cdots ({\bf e}_m\otimes 1)(1\otimes {\bf e}_m)\ . $$
Using the commutativity property  in Lemma 1.2(b), we obtain
$$\hat \Delta_{\bf i}({\bf e}_{\bf i})
=\big(({\bf e}_1\otimes 1)({\bf e}_2\otimes 1)\cdots
({\bf e}_m\otimes 1)\big)\big((1\otimes {\bf e}_1)(1\otimes {\bf
e}_2)\cdots 
(1\otimes {\bf e}_m)\big)$$
Finally, using the identities $(u\otimes 1)(v\otimes 1)=uv\otimes 1,
(1\otimes u)(1\otimes v)=1\otimes uv$ for any $u,v\in \hat {\cal U}$,
we obtain 
$\hat \Delta_{\bf i}({\bf e}_{\bf i})=({\bf e}_{\bf i}\otimes 1)
(1\otimes {\bf e}_{\bf i})={\bf e}_{\bf i}\otimes 
{\bf e}_{\bf i}$. Theorem 1.1 is proved.  $\triangleleft$

\medskip 

Now we define the {\it restricted dual} algebra 
${\cal A}={\cal U}^0$ of 
${\cal U}$. As a vector space, ${\cal A}$ is the set of all 
${\bf k}$-linear forms $x:{\cal U}\to {\bf k}$ such that  
$x$ vanishes on ${\cal U}(\gamma)$ for all but 
finitely many $\gamma\in {\bf Z}_+^r$. In other words, 
${\cal A}\cong \bigoplus_\gamma {\cal A}(\gamma)$
where ${\cal A}(\gamma)={\rm Hom}_{\bf k}({\cal U}(\gamma),{\bf k})$. 

We define the multiplication ${\cal A}\otimes {\cal A}\to {\cal A}$ 
by the formula $(xy)(u)=(x\otimes y)\big(\Delta(u)\big)$
where $(x\otimes y)\big(u_1\otimes u_2\big)=x(u_1)y(u_2)$. 
Thus, ${\cal A}$
becomes a ${\bf Z}_+^r$-graded algebra (with the unit ${\bf
k}\to {\cal A}$ dual to the counit 
$\varepsilon:{\cal U}\to {\bf k}$).

Denote by $(x,u)\mapsto x(u)$ the natural
non-degenerate evaluation pairing 
${\cal A}\times {\cal U}\to {\bf k}$. 
Furthermore, we define the pairing 
${\cal A}\times \hat {\cal U}_{\bf i}\to P_{\bf i}$ by the formula
$x(\sum t_\gamma u_\gamma)=\sum x(u_\gamma)t_\gamma.$
(The sum is finite by the definition of ${\cal A}$.) 

For every family  ${\bf E}$ as above define a map 
$\psi_{\bf i}=\psi_{\bf i,E}:
{\cal A} \to P_{\bf i}$ by the formula 
$\psi_{\bf i}(x):=x({\bf e}_{\bf i})$.
Expanding ${\bf e}_{\bf i}$ into a power series we obtain 
$$\psi_{\bf i}(x)=\sum_{a_1,\ldots,a_m\in {\bf Z}_+}
x\big(E_1^{[a_1]}E_2^{[a_2]}\cdots
E_m^{[a_m]}\big)t_1^{a_1}t_2^{a_2}\cdots t_m^{a_m}
\eqno (1.7)$$
where $E_k^{[n]}=\displaystyle{{E_k^n\over [n]_{q_k}!}}$.
Note that the sum in (1.7) is always finite 
because $x$ vanishes on all but finitely many monomials 
$E_1^{a_1}\cdots
E_{i_m}^{a_m}$. Define a ${\bf Z}_+^r$-grading on 
$P_{\bf i}$ by ${\rm deg}(t_k)=\alpha_{i_k}$ and denote 
by $P_{\bf i}(\gamma)$ 
the graded component of degree $\gamma$ in $P_{\bf i}$.   

\proclaim Corollary 1.3. For any sequence 
${\bf i}=(i_1,\ldots,i_m)$ 
and a family ${\bf E}$ of elements 
$E_k \in {\cal U}(\alpha_{i_k})$ ($k=1,\ldots,m$), the map
$\psi_{{\bf i}, {\bf E}}:{\cal A} \to P_{\bf i}$ defined by (1.7) 
is a homomorphism of ${\bf Z}_+^r$-graded algebras.

The proof of Corollary 1.3 repeats that of Corollary 0.2.
$\triangleleft$

\medskip

\noindent {\bf Remark 1}. One can prove (see e.g. [10]) 
that ${\cal A}$ 
is a ${\bf q}^t$-braided bialgebra (where ${\bf q}^t$ is the 
transpose of 
${\bf q}$). Moreover, starting with an arbitrary 
${\bf q}^t$-braided algebra
${\cal A}$, one  recovers ${\cal U}$ as the restricted 
dual of ${\cal A}$. 
So the result of Corollary 1.2 holds for any  ${\bf
q}^t$-braided bialgebra ${\cal A}$.

\noindent {\bf Remark 2}. Let 
${\cal A}_1=\bigoplus\limits_{i=1}^r {\cal
A}(\alpha_i)$.  
Corollary 1.2 implies that any morphism ${\cal A}_1 \to
\bigoplus\limits_{i=1}^r P(\alpha_i)$
of ${\bf Z}_+^r$-graded vector spaces 
extends to an algebra homomorphism. If ${\cal A}$ 
is generated by 
${\cal A}_1$, then this extension is unique. 
Thus, in the latter
case all the homomorphisms ${\cal A}\to P_{\bf
i}$ of ${\bf Z}_+^r$-graded algebras are parametrized by
the space $\bigoplus\limits_{i=1}^r \big({\cal U} 
(\alpha_i)\bigotimes P_{\bf
i}(\alpha_i)\big)$. 

\medskip 

We define the {\it universal} element ${\cal R}\in 
{\cal A}\bigotimes \hat {\cal U}$ as follows. 
The tensor product ${\cal
A}\bigotimes\hat {\cal U}$ is canonically 
identified with the space of all
linear maps ${\cal U}\to \hat {\cal U}$. 
Then ${\cal R}$ is the element in 
${\cal A}\bigotimes \hat {\cal U}$ corresponding to the inclusion
${\cal U}\hookrightarrow \hat {\cal U}$.

\proclaim Proposition 1.4. 
The element ${\cal R}$ satisfies $(\psi_{\bf i,E}\otimes 
{\rm id})({\cal R})={\bf
e}_{\bf i,E}$ 
for any ${\bf i}$ and ${\bf E}$ as above.

The proof of Proposition 1.4 coincides with that of Proposition 0.3.  
$\triangleleft$

\medskip

Clearly, the correspondence 
$$c_1t_1 + \cdots + c_mt_m \mapsto {\rm exp}_{q_1}(c_1t_1E_1)
\cdots
{\rm exp}_{q_m}(c_mt_mE_m)$$
is a map from the 
$m$-dimensional ``quantum" affine space 
$(P_{\bf i})_1=\oplus_{l=1}^m {\bf k}\cdot t_l$ 
to the set of group-like elements in $\hat {\cal U}_{\bf i}$.
This map can be regarded as a deformation of the morphism (0.1).

\medskip

Now let us turn to the fields of fractions. 
For an algebra ${\cal B}$ without zero divisors, 
${\cal F}({\cal B})$ is a vector space of right 
fractions (see Appendix).
Its elements can be written as $ab^{-1}$, 
where $a, b \in {\cal B}$ and $b \ne 0$.
As shown in the Appendix, for any sequence 
${\bf i}$ and any subalgebra 
${\cal B}\i P_{\bf i}$, the space
${\cal F}({\cal B})$ is a skew-field.
Note that $\psi_{\bf i}$ induces an embedding of 
skew fields
$\overline \psi_{\bf i}:{\cal F}({\cal A}/{\rm Ker}~ \psi_{\bf i})
\hookrightarrow {\cal F}(P_{\bf i})$. 

Now consider two elements ${\bf e}_{\bf i,E}$ and
${\bf e}_{\bf i',E'}$ corresponding via (1.6) 
to two sequences of indices 
${\bf i}=(i_1,\ldots,i_m)$ and 
${\bf i'}=(i'_1,\ldots,i'_n)$ and two families
of elements 
${\bf E}=(E_1,\ldots,E_m)$ and ${\bf E'}=(E'_1,\ldots,E'_n)$.
Let $t_1,\ldots,t_m$ (resp. $t'_1,\ldots,t'_n$) be the standard 
generators of
$P_{\bf i}$ (resp. $P_{\bf i'}$). 

\proclaim Proposition 1.5. Assume 
that ${\rm Ker}~\psi_{\bf i,E}= {\rm
Ker}~ \psi_{\bf i',E'}$, and 
$\overline \psi_{\bf i',E'}$ is an 
isomorphism of skew-fields 
${\cal F}({\cal A}/{\rm Ker}~ 
\psi_{\bf i',E'})$ and ${\cal F}(P_{\bf i'})$.  
Then the map ${\bf R}:=\overline \psi_{\bf i,E}\circ (\overline
\psi_{\bf
i',E'})^{-1}$ is an embedding ${\cal F}(P_{\bf i'}) 
\hookrightarrow {\cal F}(P_{\bf i})$,
and the following identity holds in 
${\cal F}(P_{\bf i})\bigotimes \hat
{\cal U}$:
$${\rm exp}_{q_1}(t_1E_1) \cdots 
{\rm exp}_{q_m}(t_mE_m)
={\rm exp}_{q'_1}(p_1E'_1)\cdots 
{\rm exp}_{q'_n}(p_nE'_n) \ , \eqno (1.8)$$
where $q_k=q_{i_k,i_k}$ for $k=1,\ldots,m$, 
and $q'_l=q_{i'_l,i'_l}, \,  
p_l={\bf R}(t'_l)$ for $l=1,\ldots,n$.

\noindent {\sl Proof}. We omit subscripts 
${\bf E}$ and ${\bf E}'$ in the
formulas below. Let ${\bf p}_{\bf i}:
{\cal A}\to {\cal A}/{\rm
Ker}~\psi_{\bf i}$ 
be the canonical 
projection. Denote ${\cal R}_{\bf i}
=({\bf p}_{\bf i}\otimes {\rm
id})({\cal R})$. 
Then Proposition 1.4 implies that 
$$(\overline \psi_{\bf i} \otimes 
{\rm id})({\cal R}_{\bf i})={\bf e}_{\bf
i}\eqno (1.9)$$
for every ${\bf i}\in R(w)$. We are done since 
${\bf R}\otimes {\rm id}=
(\overline  \psi_{\bf i}\otimes {\rm id})
\circ (\overline \psi_{\bf i'}\otimes {\rm id})^{-1}$. 
Proposition 1.5 is proved. $\triangleleft$

\bigskip
 
\noindent  {\bf 2. Feigin's conjecture and other results for quantum groups}

\medskip

Throughout this section we will work over the field ${\bf k}=k(q)$ where $k$ 
is a field of characteristic $0$ (say, $k={\bf C}$ as 
in the introduction), and $q$ is a variable (or a purely
transcendental element over $k$). Let $A=(a_{ij})$ be a 
symmetrizable Cartan matrix of size $r \times r$, and 
$C=(C_{ij})$ be the corresponding symmetric matrix with 
integer entries. In this section we consider a matrix 
${\bf q}$ of the form ${\bf
q}=(q_{ij})=(q^{C_{ij}})$. We denote $q_i:=q_{ii}=q^{C_{ii}}$ for all $i$. 
 
Similarly to  [8], Chapter 1, we define  
the {\it quantized enveloping algebra} ${\cal U}$ and the 
{\it quantum group} 
${\cal A}$ associated with $A$ as follows. 
First, let $\overline {\cal U}$ be the free algebra over 
$k(q)$ generated by $E_1,\ldots,E_r$. We make 
$\overline {\cal U}$ into a ${\bf q}$-braided 
bialgebra (see Section 1).  
Second, the restricted dual algebra 
$\overline {\cal A}$ of $\overline {\cal U}$ is defined 
as in Section 1. 
Define a homomorphism $f:\overline {\cal U}\to \overline 
{\cal A}$ by $f(E_i)= x_i$ where $x_i$ is the only 
element in $\overline {\cal A}(\alpha_i)$ 
such that $x_i(E_i)=1$. Finally, define 
${\cal U}:=\overline {\cal U}/{\rm Ker}~f$ and 
${\cal A}:={\rm Im}~f$, and keep the above notation 
for the generators.  In particular, ${\cal U}\cong 
{\cal A}$ via $E_i\mapsto x_i$. 
It is well-known that the right kernel of the evaluation 
pairing ${\cal A}\otimes \overline {\cal U}\to k(q)$
coincides with ${\rm Ker}~f$. Hence the induced pairing 
$${\cal A}\otimes {\cal U}\to k(q) \eqno (2.1)$$
is non-degenerate, so we identify ${\cal A}$ 
with the restricted dual 
algebra to the ${\bf q}$-braided bialgebra ${\cal U}$ 
(and denote the 
evaluation pairing (2.1) by $(x,E)\mapsto x(E)$). 
Note that the generators $E_1,\ldots,E_r$ of ${\cal U}$ 
(as well as the generators $x_1,\ldots,x_r$ of ${\cal A}$) 
are subject to  
the quantum Serre 
relations ([8], Section 1.4.3, or Section 3 below).

The algebra ${\cal U}$ is ${\bf Z}_+^r$-graded via 
${\rm deg}~ E_i=\alpha_i$. 
The pairing ${\cal A}\times {\cal U}\to k(q)$ extends to 
the $P_{\bf i}$-linear pairing ${\cal A}\times \hat 
{\cal U}_{\bf i}\to k(q)$ 
(we denote it  by $(x,u)\mapsto x(u)$), where $\hat 
{\cal U}_{\bf i}:=P_{\bf
i}\bigotimes \hat {\cal U}$ and  $P_{\bf i}$ is a 
$k(q)$-algebra generated by
$t_1,\ldots,t_m$ subject to the relations 
$t_lt_k=q^{C_{i_k,i_l}}t_kt_l$
for $1\le k<l\le m$.
 
For the convenience of the reader, we summarize the results
from Section 1 for the quantum groups ${\cal A}$ and ${\cal U}$ in
the following theorem.  
 
\proclaim Theorem 2.1. Let ${\bf i}=(i_1,\ldots,i_m)$ be 
any sequence. Then \item {(a)} the element   
$${\bf e}_{\bf i}={\rm exp}_{q_{i_1}}(t_1E_{i_1}) \cdots
{\rm exp}_{q_{i_m}}(t_mE_{i_m})$$ 
is a group-like element in $\hat {\cal U}_{\bf i}$; 
\item {(b)} the element ${\bf e}_{\bf i}$ gives rise to an algebra 
homomorphism 
$\psi_{\bf i}:{\cal A}\to P_{\bf i}$ defined by 
$\psi_{\bf i}(x):=x({\bf e}_{\bf i})$; 
\item {(c)} there is a unique element 
${\cal R}\in {\cal A}\bigotimes
\hat {\cal U}$ satisfying 
$(\psi_{\bf i}\otimes {\rm id})({\cal R})=
{\bf e}_{\bf i}$ for all ${\bf i}$;
\item {(d)} the homomorphism $\psi_{\bf i}$ satisfies 
$$\psi_{\bf i}(x)=
\sum_{a_1,\ldots,a_m\ge 0} x(E_{i_1}^{[a_1]}E_{i_2}^{[a_2]}\cdots
E_{i_m}^{[a_m]})t_1^{a_1}t_2^{a_2}\cdots t_m^{a_m} \ , \eqno (2.2)$$
where $\displaystyle{E_i^{[n]}={1\over
[n]_{q_i}!}E_i^n}$, and $[n]_q!$ is defined in (0.4). 
In particular, for $i=1,\ldots,r$ we 
have $\psi_{\bf i}(x_i)=\sum\limits_{k:i_k=i} t_k$, and this 
determines $\psi_{\bf i}$ uniquely; 
\item {(e)} ${\rm Ker}~~\psi_{\bf i}=\{x\in {\cal A}:
x(E_{i_1}^{a_1}E_{i_2}^{a_2}\cdots E_{i_m}^{a_m})=0~{\rm for
~all}~a_1,\ldots,a_m\in {\bf Z}_+\}.$

{\bf Remark}. After the identification ${\cal A}\cong {\cal U}$ 
as above,  $\psi_{\bf i}$ coincides with Feigin's 
homomorphism 
$\Phi({\bf i}):{\cal U}\to P_{\bf i}$. B. Feigin introduced this 
homomorphism in his 
talk at RIMS in 1992 (see e.g. [5] and [6]). 
  
\vskip .2cm

Let $W$ be the Weyl group associated with the Cartan matrix $A$. By
definition, $W$ is generated by simple reflections
$s_1,\ldots, s_r:{\bf Z}^r\to {\bf Z}^r$ where
$s_i(\alpha_j)=\alpha_j-a_{ij}\alpha_i$. We call a 
sequence ${\bf i}=(i_1,\ldots,i_m)$ of indices a 
{\it reduced expression} of $w\in W$ if 
$w=s_{i_1}s_{i_2}\cdots s_{i_m}$,
and the above expression of $w$ is the shortest 
(we call ${\bf i}$ simply a
reduced expression if $w$ is not specified). 
We set $l(w):=m$ and call
$l(w)$ the {\it length} of $w$. Denote by $R(w)$ 
the set of all reduced
expressions  of $w$. It is well-known that 
$W$ is a Coxeter group, so the defining 
relations between $s_1,\ldots,s_r$ are of the form 
$(s_is_j)^l=1$ where 
$l\in \{2,3,4,6\}$. (More precisely, 
$l=a_{ij}a_{ji}+2$ if $a_{ij}a_{ji}<3$ 
and $l=6$ if $a_{ij}a_{ji}=3$.) 
It follows that every two reduced expressions  
of an element $w\in W$ are connected by a chain of moves 
$$({\bf i}_1,(i,j,i,\ldots),{\bf i}_2)\mapsto 
({\bf i}_1,(j,i,j\ldots),{\bf i}_2)$$
where each fragment in parentheses has the length $l$. 
If the Weyl group $W$ is finite then there is a unique 
element of the maximal length in $W$ which we 
denote by $w_0$.  

Let us study the kernel of $\psi_{\bf i}$. 
According to Theorem 2.1(e), ${\rm Ker}~\psi_{\bf i}$ 
is the orthogonal complement of 
the subspace ${\cal U}({\bf i})\i {\cal U}$ spanned by all 
monomials $E_{i_1}^{a_1}E_{i_2}^{a_2}\cdots
E_{i_m}^{a_m}$.

\proclaim Lemma 2.2. \item {(a)} For every sequence ${\bf i}$ 
there is a reduced expression ${\bf i}'$ such that 
${\cal U}({\bf i})={\cal U}({\bf i}')$.
 Moreover, ${\bf i'}$ can always be chosen as a
subsequence of ${\bf i}$.  
\item {(b)}  For any $w\in W$ and ${\bf i},{\bf i}'\in R(w)$, 
 we have ${\cal U}({\bf i})={\cal U}({\bf i}')$. 
\item {(c)} ${\cal
U}({\bf i})$ contains the subalgebra in ${\cal U}$ generated by all 
$E_i$ such that $l(ws_i)=l(w)-1$. Therefore, 
${\cal U}({\bf i})={\cal U}$
for every ${\bf i}\in R(w_0)$ when $W$ is finite. 

\vskip .1cm

\noindent {\sl Proof}. The collection of the subspaces 
$\{{\cal U}({\bf i})\}$ is a 
multiplicative semigroup with respect to the product of vector
subspaces in ${\cal U}$. By definition, 
$${\cal U}(i_1,i_2,\ldots,i_m)=
{\cal U}(i_1)\cdot {\cal U}(i_2)\cdots 
{\cal U}(i_m)$$
where ${\cal U}(i)$ is a subalgebra in ${\cal U}$ generated 
by $E_i$, 
$i=1,\ldots,r$.

We have ${\cal U}(i){\cal U}(i)={\cal U}(i)$, and for every 
pair $(i,j)$ with
$a_{ij}a_{ji}<4$ the following relation holds:
$${\cal U}(i)\cdot {\cal U}(j)\cdot {\cal U}(i)\cdots = {\cal U}(j)\cdot
{\cal U}(i)\cdot {\cal U}(j)\cdots\eqno (2.3)$$
where each product contains $l$ factors. The identity (2.3) can be 
proved by the standard arguments for the
algebras ${\cal U}$ whose Cartan matrices are of types 
$A_1\times A_1,A_2,B_2$ or $G_2$.  
 
It follows that every ${\cal U}({\bf i})$ equals to  
${\cal U}({\bf i'})$ for some reduced subsequence 
${\bf i}'$ of ${\bf i}$, which 
proves (a). Part (b) also follows because the braid
relations (2.3) can be used to move from any 
reduced expression of $w\in W$ to any other one. 

(c) Let $J=J_w$ be the the set of all $i$ satisfying $l(ws_i)=l(w)-1$. 
For each $i\in J$, there exists ${\bf i}\in R(w)$ such that 
${\bf i}$ ends with $i$. Using (b) we see that ${\cal U}({\bf i})E_i\i 
{\cal U}({\bf i})$ for any ${\bf i}\in R(w)$, $i\in J$. This completes 
the proof of  Lemma 2.2. $\triangleleft$

\medskip 

We define $I_w:={\rm Ker}~\psi_{\bf i}$ 
for any ${\bf i}\in R(w)$. Since ${\cal A}/I_w$ is 
isomorphic to $\psi_{\bf i}({\cal A})$, it follows that 
${\cal F}({\cal A}/I_w)$ is a skew field (see Appendix).

\proclaim Theorem 2.3. For every  $w\in W$ and 
${\bf i}\in R(w)$ the map $\psi_{\bf i}$ induces an
isomorphism of skew fields 
$$\overline \psi_{\bf i}:{\cal F}({\cal A}/I_w)\cong {\cal F}
(P_{\bf i}) \ .  \eqno (2.4)$$

Taking $w=w_0$ we obtain the following 

\proclaim Corollary 2.4. (Feigin's conjecture). For any ${\bf
i}\in R(w_0)$, the homomorphism  $\psi_{\bf i}:{\cal A}\to P_{\bf i}$ 
is an embedding, and it induces an isomorphism
of skew-fields
$$\overline \psi_{\bf i}:{\cal F}({\cal A})\cong 
{\cal F}(P_{\bf i}) \ .$$

\vskip .2cm
\noindent {\sl Proof of Theorem 2.3}.  It is enough to prove that 
for any ${\bf i}\in R(w)$ the image of $\psi_{\bf i}$ generates 
${\cal F}(P_{\bf i})$, that is, $t_1,\ldots,t_m$ belong to 
${\cal F}({\rm Im}~ \psi_{\bf i})$.
We will deduce  this statement from Proposition 0.6.  
Then we neeed the following.
 
\proclaim Proposition 2.5. 
For each element $x\in {\cal A}$ satisfying 
$$\psi_{\bf i}(x)=t_1^{a_1}t_2^{a_2}\cdots t_m^{a_m}\eqno (2.5)$$ 
with $a_1>0$, there is an element $y\in
{\cal A}$ such that $\psi_{\bf i}(y)=ct_1^{a_1-1}t_2^{a_2}
\cdots t_m^{a_m}$ where $c\in k(q),c\ne 0$.

\noindent {\sl Proof of Proposition 2.5}. 
For $i=1,\ldots,r$, let $E_i^*:{\cal A}\to {\cal A}$ be 
the adjoint operator of the left multiplication operator 
$E\mapsto E_iE$ in 
${\cal U}$. Thus, 
the element $E_i^*(x)$  is
determined by 
the equations  $(E_i^*(x))(E)=x(E_iE)$ 
for every $E\in {\cal U}$.

We will show that $y$ can be chosen as $y=E_i^*(x)$.
    
Indeed, (2.5) means that the right hand side of the 
expansion (2.2) for $\psi_{\bf i}(x)$ reduces to one summand or, 
equivalently, 
$$x(E_{i_1}^{[b_1]}\cdots E_{i_m}^{[b_m]})=0\eqno
(2.6)$$
unless $(b_1,\ldots,b_m)=(a_1,\ldots,a_m)$.

By (2.2), we have   
$$\psi_{\bf i}(y)=\psi_{\bf i}(E^*_{i_1}(x))
=\sum_{b_1,\ldots,b_m\in {\bf Z}_+}  
(E_{i_1}^*(x))\big(E_{i_1}^{[b_1]}E_{i_2}^{[b_2]}\cdots
E_{i_m}^{[b_m]}\big)t_1^{b_1}t_2^{b_2}\cdots t_m^{b_m}$$
$$=\sum_{b_1,\ldots,b_m\in {\bf Z}_+}
x\big(E_{i_1}E_{i_1}^{[b_1]}E_{i_2}^{[b_2]}\cdots
E_{i_m}^{[b_m]}\big)t_1^{b_1}t_2^{b_2}\cdots t_m^{b_m} \ .$$

In view of (2.6), 
$$\psi_{\bf i}(y)=x(E_{i_1}E_{i_1}^{[a_1-1]}
E_{i_2}^{[a_2]}\cdots E_{i_m}^{[a_m]}) t_1^{a_1-1}t_2^{a_2}\cdots
t_m^{a_m}=ct_1^{a_1-1}t_2^{a_2}\cdots t_m^{a_m}$$ 
with $c\ne 0$ as desired.  $\triangleleft$

\medskip

Taking $x$ and $y$ as in Proposition 2.5, we see that
$$t_1=c\psi_{\bf i}(x)(\psi_{\bf i}
(y))^{-1}\in ~{\cal F}({\rm Im}~\psi_{\bf i}) \ . \eqno (2.7)$$

To complete the proof of Theorem 2.3, we proceed 
by induction on $m$. If $m=1$ then 
$t_1\in {\rm Im}~\psi_{\bf i}=P_{\bf i}$. 
So let $m\ge 2$, denote ${\bf i}'=(i_2,\ldots,i_m)$ and assume 
that Theorem 2.3 holds for ${\bf i}'$, that is,     
$$t_2,t_3,\ldots,t_m\in {\cal F}({\rm Im}~\psi_{\bf i'}) \ . 
\eqno (2.8)$$

Note that $P_{{\bf i}'}$ is naturally embedded into 
$P_{\bf i}$ as a subalgebra generated by $t_2,\ldots,t_m$. 
 
In view of the formula for $\psi_{\bf i}(x_i)$ in Theorem 2.1(d),
$$\psi_{{\bf i}'}(x_i)=
\psi_{\bf i}(x_i)~ (i\ne i_1),~\psi_{{\bf i}'}(x_{i_1})=
\psi_{\bf i}(x_{i_1})-t_1 \ . \eqno (2.9)$$

Using (2.7). we see that  
$$\psi_{{\bf i}'}(x_i)\in {\cal F}({\rm Im}~\psi_{\bf i}) \ ,
(i=1,\ldots,r)$$ 
hence   
$${\cal F}({\rm Im}~\psi_{{\bf i}'})\i  
{\cal F}({\rm Im}~\psi_{\bf i}) \ .$$

Combining this with the inductive assumption (2.8), we conclude that 
$t_2,\ldots,t_m\in {\cal F}({\rm Im}~\psi_{\bf i})$.
Since $t_1$ also belongs to ${\cal F}({\rm Im}~\psi_{\bf i})$, 
Theorem 2.3 is proved.  $\triangleleft$

\medskip 
 
\noindent {\bf Proof of Proposition 0.8.} 
let ${\cal B}$ be the algebra of the upper triangular 
$n\times n$-matrices 
over ${\bf C}(q)$ (with the unity $I$, the identity matrix). Let 
$\rho:{\cal U}\to {\cal B}$ be a 
representation of ${\cal U}$ given by
$\rho(E_i)=E_{i,i+1}$,
where $E_{ij}$ is the matrix unit. The representation 
$\rho$ extends naturally to  ${\rm id}\otimes 
\rho:{\cal A}\bigotimes {\cal U}\to  {\cal A}\bigotimes {\cal B}$. 
We identify the latter algebra with ${\cal B}({\cal A})$, 
the algebra of upper triangular matrices over ${\cal A}$.  

\proclaim Lemma 2.6. We have $X=({\rm id}\otimes 
\rho) ({\cal R})$, where ${\cal R}$ is the universal 
element in ${\cal A}\bigotimes \hat {\cal U}$. 

\noindent {\sl Proof}. Note that ${\cal B}$ is a 
${\bf Z}_+^{n-1}$-graded algebra via ${\rm deg} (E_{ii})=0,~{\rm deg}
(E_{ij})=\alpha_{ij}=\alpha_i+\cdots \alpha_{j-1}$ 
($1\le i<j\le n$),  
and $\rho$ preserves the ${\bf Z}_+^{n-1}$-grading ($n=r+1$).
Therefore,  the formula (0.9) for ${\cal R}$ implies
$$({\rm id}\otimes 
\rho) ({\cal R})=I+\sum\limits_{i<j} 
\sum\limits_{b\in B_{ij}} b^*\rho(b) \ . $$ 
where $B_{ij}$ is a basis in ${\cal U}(\alpha_{ij})$ 
and $\{b^*\}$ is the dual basis in ${\cal A}(\alpha_{ij})$.  
We choose $B_{ij}$ to consist of the products (in any order) 
of the generators
$E_i,E_{i+1},\ldots, E_{j-1}$. It is easy to see that 
$\rho(b)=E_{ij}$ for the element 
$b=b_{ij}=E_iE_{i+1}\cdots E_{j-1}$ in $B_{ij}$, and 
$\rho(b)=0$ if $b\in B_{ij}, b\ne b_{ij}$.  
 Denote $x_{ij}=(b_{ij})^*$. To identify these $x_{ij}$ with 
those in Section 0 we have to verify the relations (0.11). 
As an algebra, ${\cal A}$ is generated by the 
$x_i:=x_{i,i+1}$ ($i=1,\ldots,r$) 
subject to the quantum Serre relations. The relations (0.11) can 
be verified similarly 
to those  between the $t_{ij}$ in [2], Section 3 
(they also follow from the relations in [4]).  Thus,  
$({\rm id}\otimes 
\rho) ({\cal R})=I+\sum\limits_{i<j}
x_{ij}E_{ij}=X$. Lemma 2.6 is proved. $\triangleleft$ 

\medskip 

To complete the proof of Proposition 0.8, note that, 
for every ${\bf i}$ we have     
$\psi_{\bf i}(X)=(\psi_{\bf i}\otimes 
{\rm id})\circ ({\rm id}\otimes \rho)({\cal R})
=({\rm id}\otimes \rho)\circ (\psi_{\bf i}\otimes 
{\rm id})({\cal R})=({\rm id}\otimes \rho)({\bf e}_{\bf i})$, 
by (0.8). 
The formula (0.12) follows since $({\rm id}\otimes \rho)\big 
({\rm exp}_{q_{i_k}}(t_kE_{i_k})\big )=I+t_kE_{i_k,i_k+1}$ 
for all $k$. Proposition 0.8 is proved. $\triangleleft$ 

\medskip
 
We have the following obvious corollary of Theorem 2.3.

\proclaim Corollary 2.7. 
For any $w\in W$ and ${\bf i}$,${\bf i}'\in R(w)$ 
there is an isomorphism of skew-fields
$$R_{\bf i}^{\bf i'}: 
{\cal F}(P_{\bf i})~\widetilde \to ~ {\cal F}(P_{\bf i'}) 
\eqno (2.10)$$
defined by 
$R_{\bf i}^{\bf i'}:=\overline 
\psi_{\bf i'}\circ (\overline \psi_{\bf i})^{-1}$.

The ``transition 
maps" $R_{\bf i}^{\bf i'}$ lead to 
identities between quantum exponentials given by Corollary 0.10. 
To compute each $p_k$ in (0.14), it 
is enough (in principle) to do this for the following pairs:
$${\bf i}=(i,j,i\ldots), ~{\bf i'}=(j,i,j,\ldots) \ ,$$ 
of the length $l$ each, where $l$ is the order of $s_is_j$ in $W$.  
Recall that $l=2,3,4$, or $6$ for any Weyl group. 
In the following proposition, we 
compute ${\bf R}_{\bf i'}^{\bf i}$ for these ${\bf i},{\bf i'}$ 
with $l=2,3,$ or $4$ 
(when $l=6$ the explicit expressions for $p_k$ are more 
complicated, so we do not present them here).

\proclaim Proposition 2.8. Let $t_1,\ldots,t_l$ 
(resp. $t'_1,\ldots,t'_l$) be 
standard generators of $P_{(iji\ldots)}$  
(resp. $P_{(jij\ldots)}$). 
We denote $p_k={\bf R}_{(jij\cdots )}^{(iji\cdots )}(t'_k)$ 
$(k=1,\ldots,l$).  
\item {(a)} If $l=2$ then $(p_1,p_2)= (t_2,t_1)$. 
\item {(b)} If $l=3$ then $(p_1,p_2,p_3)=
\big(t_2t_3(t_1+t_3)^{-1},~t_1+t_3,~(t_1+t_3)^{-1}t_1t_2\big)$.   
\item {(c)} If $l=4$ and $a_{ij}=-2,a_{ji}=-1$ then 
$p_1,p_2,p_3,p_4$ are determined by the following equations:
$$p_2p_3=t_1t_2+t_1t_4+t_3t_4,~~p_2p_3p_4=t_1t_2t_3 \ , $$
$$p_2^2p_3=t_1^2t_2+(t_1+t_3)^2t_4,~~p_1p_2^2p_3=t_2t_3^2t_4 \ .$$ 

{\noindent {\sl Proof}. (a) By definition, 
$\psi_{(ij)}:(x_i,x_j)\mapsto (t_1,t_2)$ and 
$\psi_{(ji)}:(x_i,x_j)\mapsto (t'_2,t'_1)$. Thus, 
$p_1=t_2$, and $p_2=t_1$ as claimed.  

Part (b) is proved in Section~0 (see Example 1).

(c)  Let ${\cal U}_{ij}$ be 
the subalgebra of ${\cal U}$ generated 
by $E_i$ and $E_j$. 
Define ${\cal B}_{ij}$ as the quotient algebra of ${\cal U}_{ij}$ 
modulo the relations $E_i^3=E_j^2=0$ 
(we keep the same notation for generators). 
It is easy to see ([8], or Section 3 below) that the following are
all the defining relations in ${\cal B}_{ij}$: 
$$E_i^3=E_j^2=E_jE_iE_j=0,~E_i^2E_jE_i=E_iE_jE_i^2 \ .$$
Using these relations, it is easy to prove that the homogeneous
components ${\cal B}_{ij}(\alpha_i+\alpha_j)$, 
${\cal B}_{ij}(2\alpha_i+\alpha_j)$, and 
${\cal B}_{ij}(2\alpha_i+2\alpha_j)$ of ${\cal B}_{ij}$
have the following bases: $\{E_iE_j,~E_jE_i\}$ for 
${\cal B}_{ij}(\alpha_i+\alpha_j)$, $\{E_i^2E_j,~E_iE_jE_i,~E_jE_i^2\}$ for 
${\cal B}_{ij}(2\alpha_i+\alpha_j)$, and 
$\{E_jE_i^2E_j\}$ for 
${\cal B}_{ij}(2\alpha_i+2\alpha_j)$.

Applying the projection $\rho:{\cal U}_{ij}\to {\cal B}_{ij}$ to 
both sides of (0.14), we obtain the following  equation in 
${\cal F}(P_{ijij})\bigotimes {\cal B}_{ij}$: 
$$(1+p_1E_j)
(1+p_2E_i+{p_2^2E_i^2\over 1+q^2})(1+p_3E_j)
(1+p_4E_i+{p_4^2E_i^2\over 1+q^2})$$
$$=(1+t_1E_i+{t_1^2E_i^2\over 1+q^2})(1+t_2E_j)
(1+t_3E_i+{t_3^2E_i^2\over 1+q^2})(1+t_4E_j) \ . \eqno (2.11)$$
The desired expressions for 
$p_2p_3, ~p_2p_3p_4,~p_2^2p_3$, and $p_1p_2^2p_3$ can be obtained 
from (2.11) by comparing the coefficients of 
$E_iE_j, ~E_iE_jE_i,~E_i^2E_j$, 
and $E_jE_i^2E_j$  
respectively on both sides of (2.11). Proposition 2.9 is
proved. $\triangleleft$

\vskip .2cm

{\bf Remark}. Taking in the identities of Proposition 2.8 the homogeneous  components of degrees $\alpha_i+(1-a_{ij})\alpha_j$ and
$\alpha_j+(1-a_{ji})\alpha_i$ yields quantum Serre relations 
between $E_i$ and $E_j$.

\vskip .3cm 

We conclude this section by a proof of Proposition 0.11. 
For a skew-symmetric $m\times m$-matrix $S=(S_{kl})$ with 
integer entries let $P_S$ be a $k(q)$-algebra generated by
$t_1,\ldots,t_m$ subject to the relations
$$t_lt_k=q^{S_{kl}}t_kt_l \ . \eqno (2.12)$$ Note that 
$P_{\bf i}=P_{S({\bf i})}$ where the matrix 
$S({\bf i})$ is defined in (0.16). Recall that two 
$m\times m$-matrices $S$ and $S'$ are called equivalent 
if there is a matrix $T=T_{kl}\in SL_m({\bf Z})$ such that $S'=TS T^t$. 
It is easy to see that ${\cal F}(P_S)\cong {\cal F}(P_{S'})$ 
if $S$ and $S'$ are equivalent: one can choose such an isomorphism 
${\cal F}(P_{S'})\widetilde 
\to  {\cal F}(P_S)$ by sending each generator $t'_k$ of $P_{S'}$ to 
the monomial $t_1^{T_{k,1}}t_2^{T_{k,2}}\cdots t_m^{T_{k,m}}$ 
in the generators of $P_S$. The converse statement was proved by A. Panov.

\proclaim Proposition 2.9 ([11], Theorem 2.19). 
Let $S,S'$ be skew-symmetric $m\times m$ matrices 
with the integer entries. 
Then ${\cal F}(P_S)\cong {\cal F}(P_{S'})$ if and only 
if $S$ and $S'$ are equivalent.

Thus, Proposition 2.9 means that the existence of {\it any} isomorphism 
${\bf R}:{\cal F}(P_{S'})\to {\cal F}(P_S)$ implies 
that of a {\it monomial} isomorphism ${\bf M}:{\cal F}(P_{S'})
\to {\cal F}(P_S)$, that is, ${\bf M}$ 
takes each generator $t'_k$ of $P_{S'}$ to a monomial in generators 
$t_1,\ldots,t_m$ of $P_S$. 
Taking $S=S({\bf i}), S'=S({\bf i'})$, 
and ${\bf R}={\bf R}_{\bf i'}^{\bf i}:{\cal F}(P_{\bf i'})
\to {\cal F}(P_{\bf i})$ with ${\bf i},{\bf i'}\in R(w)$ for 
some $w\in W$, we obtain, 
in particular, the statement of Proposition 0.11.
Note that ${\rm R}_{\bf i'}^{\bf i}$  
is not monomial in general.

One can prove that there exists a {\it local} monomial isomorphism ${\bf
M}={\bf M}_{\bf i,i'}:{\cal F}(P_{\bf i'}) \to {\cal F}(P_{\bf
i})$. Namely, for ${\bf i},{\bf i'}$ of the form 
$${\bf i}=(i_1,\ldots,i_{a-l};i,j,i,\ldots;i_{a+1},\ldots,i_m)$$
$${\bf i'}=(i_1,\ldots,i_{a-l};j,i,j,\ldots;i_{a+1},\ldots,i_m) \ ,$$
${\bf M}(t'_k)=t_k$ if $k\le a-l$ or $k>a$, and each ${\bf M}(t'_k)$ for 
$k=a-l+1,\ldots,a$ depends only on $t_{a-l+1},\ldots,t_a$.  We will 
present such ${\bf M}$ elswhere. 

\smallskip

By Proposition 2.9, the equivalence class of $S({\bf i})$ 
for ${\bf i}\in R(w)$  depends only on $w$. If we choose some
representative ${\bf S}(w)$ of this class then, by Theorem 2.3, there is an
isomorphism  
$${\cal F}\big({\cal A}/I_w\big )\widetilde \to 
{\cal F}(P_{{\bf S}(w)}) \ . \eqno (2.16) $$

The well-known normal form for skew-symmetric matrices shows that 
${\bf S}(w)=({\bf S}_{kl})$ can be chosen uniquely  subject to the 
following requirements:
\item {(i)}  ${\bf S}_{kl}=0$ unless $k+l\ne m+1$; 
\item {(ii)}  there is a sequence $c_1,c_2,\ldots$ of 
nonnegative integers such that    
${\bf S}_{k,m+1-k}=c_1c_2\cdots c_k$
for $1\le k\le {m\over 2}$. 

\smallskip 
It would be interesting to compute the invariants $c_1,\ldots,c_m$ in 
terms of $w$, and to find a direct way to describe the isomorphism (2.16).

\bigskip 
 
\noindent {\bf 3. Extremal vectors in ${\cal A}$ and proof of 
Proposition 0.6}

\medskip

We retain terminology and notation of Section 2.
Recall that $\alpha_1,\ldots,\alpha_r$ is the standard basis in 
${\bf Z}_+^r$. 
We define a bilinear form in ${\bf Z}^r$ by the formula
$(\alpha_i,\alpha_j)=C_{ij}$ for all $i,j$.

Let us fix $\lambda=(l_1,\ldots,l_r)\in {\bf Z}_+^r$. 
For $i=1,\ldots,r$  
define linear operators $F_i=F_{i,\lambda}: 
{\cal A}\to {\cal A}$ by the formula:
$$F_i\cdot x= {v_i^{l_i}q^{-(\gamma,\alpha_i)}x x_i-v_i^{-l_i}x_ix
\over v_i-v_i^{-1}}  \eqno (3.1)$$
for $x \in {\cal A}(\gamma)$ (where $v_i=q^{C_{ii}\over 2}$).

We identify $\lambda$ with a linear form on the {\it coroot}
lattice ${\bf Z}\alpha_1^\vee\oplus\cdots 
\oplus {\bf Z}\alpha_r^\vee$ 
defined by $\lambda(\alpha_i^\vee):=l_i$ (recall that $\alpha_i^\vee=
{2\alpha_i\over (\alpha_i,\alpha_i)}$). For each reduced ${\bf
i}=(i_1,\ldots,i_m)$ we define a 
sequence of integers 
$a_1,\ldots, a_m$ by the formula
$$a_1=\lambda(\alpha_{i_1}^\vee),~~a_2=
\lambda(s_{i_1}(\alpha_{i_2}^\vee)),\ldots,
a_m=\lambda\big(s_{i_1}s_{i_{2}}\cdots 
s_{i_{m-1}}(\alpha_{i_m}^\vee)\big) \ . \eqno (3.2)$$
It is well-known that $a_k\in {\bf Z}_+$ for all $k$.

Define the  element $v({\bf i})=v({\bf i})^\lambda\in {\cal A}$
by:
$$v({\bf i})=
F_{i_m}^{a_m}F_{i_{m-1}}^{a_{m-1}}\cdots 
F_{i_1}^{a_1}\cdot 1 \ .\eqno (3.3)$$
 
The  following result refines  Proposition 0.6.

\proclaim Theorem 3.1.  In the above notation, we have 
$\psi_{\bf i}(v({\bf i}))=ct_1^{a_1}t_2^{a_2}\cdots 
t_m^{a_m}$ where $c\in k(q),c\ne 0$.

Proposition 0.6 follows by taking any $\lambda$ with
$\lambda(\alpha_{i_1}^\vee)=l_{i_1}>0$ (and $x:=c^{-1}v({\bf i})$).
 
\smallskip

\noindent {\sl Proof of Theorem 3.1}.
Let us reformulate our statement 
in terms of modules over the quantized enveloping algebra 
${\bf U}=U_q({\bf g})$. 
The $k(q)$-algebra ${\bf U}$ is generated by
$F_1,\ldots,F_r$, $E_1,\ldots,E_r$ and the invertible pairwise commuting
elements $K_1,\ldots,K_r$  subject to the following relations (see [8]):
$$K_iE_jK_i^{-1}=q^{C_{ij}}E_j,
~K_iF_jK_i^{-1}=q^{-C_{ij}}E_jK_i,~~~
E_iF_j-F_jE_i=\delta_{ij}\displaystyle{K-K^{-1} \over
v_i-v_i^{-1}} \ ;\eqno (3.4)$$ 
$$\sum_{p+p'=1-a_{ij}} (-1)^p v_i^{pp'}  
E_i^{[p]}E_j E_i^{[p']}=0, \,\, \sum_{p+p'=1-a_{ij}} (-1)^p v_i^{pp'}  
F_i^{[p]}F_j F_i^{[p']}=0 \ . \eqno (3.5)$$
The relations (3.5) are {\it quantum  Serre relations};
they hold for all $i\ne j$ where $E_i^{[n]}$ means the same as in (2.2),  
and $v_i=q^{C_{ii}\over 2}$.   

The  between $E_1,\ldots,E_r$ 
(the same relations between $F_1,\ldots,F_r$):

We identify the subalgebra generated by 
$E_1,\ldots,E_r$ with ${\cal U}$. Note that ${\bf U}$ is a 
${\bf Z}^r$-graded algebra  via 
${\rm  deg}(K_i)={\rm  deg}(K_i^{-1})=0,~~
{\rm  deg}(E_i)=\alpha_i,~{\rm  deg}(F_i)=-\alpha_i$ for 
$i=1,\ldots,r$. Note also, 
that each triple $(E_i,F_i,K_i)$ generates the
subalgebra in ${\bf U}$ isomorphic to $U_{v_i}(sl_2)$.

We will denote by the same symbol $E_i$ the operator 
$E_i:{\cal A}\to {\cal
A}$ adjoint of  the operator of the right multiplication 
$E\to EE_i$ in ${\cal U}$; for every $x\in {\cal A}$ 
the element $E_i\cdot x\in {\cal A}$ is defined  by the equations
$(E_i\cdot x)(E)=x(EE_i)$ for all $E\in {\cal U}$. 
We also define the operator $K_i:{\cal A}\to {\cal A}$ 
(depending on $\lambda$) by 
$$K_i\cdot x=K_{i,\lambda}\cdot x:= q^{\lambda(\alpha_i)
-(\gamma,\alpha_i)}x \  \eqno (3.6) $$
for all $x\in {\cal A}(\gamma)$, $i=1,\ldots,r$. 

The following result is well-known, 
(for type $A_r$ it can be found e.g, in [3]). 

\proclaim Proposition 3.2. For every $\lambda$ as above, the 
operators $F_i$ defined in (3.1) together with the 
$K_i$ and $E_i$, give rise to an action 
${\bf U}\times {\cal A}\to {\cal A}$. 
 
Denote by $V_\lambda$ the cyclic ${\bf U}$-submodule in 
${\cal A}$ 
generated by the unit $1\in {\cal A}$. It is well-known (cf. [3], [8])
that $V_\lambda$ is an integrable simple ${\bf U}$-module. 
Note also that the vector $v=1\in V_\lambda$ 
is a {\it highest weight vector of weight} 
$\lambda$ since  $E_i\cdot v=0$ and $K_i(v)=
q^{\lambda(\alpha_i)}v$ for all $i$.  
 
It is also known (see [8]), that the element $v({\bf i})$ given by 
(3.3) depends only on $w$. 
Such elements are called {\it extremal vectors} in $V_\lambda$. 
Denote ${\bf i}_k:=(i_1,\ldots,i_k)$ for $k=1,\ldots,m$. It is also 
well-known that for all $k$ we have 
$$F_{i_k}\cdot v({\bf i}_k)=0,~ E_{i_k}^a\cdot v({\bf i}_k)=
c_{k;a}F_{i_k}^{a_k-a}v({\bf i}_{k-1}),~ E_{i_k}\cdot v({\bf i}_{k-1})=0
\eqno (3.7)$$
for some $c_{k;a}\in k(q) \setminus \{0\}$, and $a=0,1,\ldots,a_k$ 
(with the agreement  $v({\bf i}_0):=v$). 
 
In view of (2.2), Theorem 3.1 is equivalent to the following.

\proclaim Proposition 3.3. There is a unique sequence 
$b=(b_1,\ldots,b_m)$ such that 
$$E_{i_1}^{b_1}\cdots E_{i_m}^{b_m}
\cdot v({\bf i})=cv\eqno (3.8)$$  
for some nonzero scalar $c\in k(q)$, namely, 
$(b_1,\ldots,b_m)=(a_1,\ldots,a_m)$, where $a_1,\ldots,a_m$ are 
given by (3.2).

\noindent {\sl Proof of Proposition 3.3}. We proceed by
induction on $m$. 

Assume that our statement is true for every reduced expression 
of length $<m$, in particular, for ${\bf i}_{m-1}$.

\smallskip 

\noindent {\bf Step 1}. Let us prove that the equality (3.8) implies that 
$b_k=a_k$ for all $k$ with $i_k\ne i_m$.

We will use the following identity in ${\bf U}$, which is 
a straightforward consequence of the relations (3.4):
$$E_{i_1}^{b_1}E_{i_2}^{b_2}\cdots E_{i_m}^{b_m}F_{i_m}^{a_m}=
\sum_{b'} F_{i_m}^{|b'|-|b|+a_m}
E_{i_1}^{b'_1}E_{i_2}^{b'_2}\cdots E_{i_m}^{b'_m}p_{b'}
\eqno (3.9)$$
where the sum is over all $b'=(b'_1,\ldots,b'_m)\in {\bf Z}_+^m$ 
such that $b'_k=b_k$ if
$i_k\ne i_m$, and $b'_k\le b_k$ if $i_k=i_m$; each $p_{b'}$ is a 
Laurent polynomial of $K_{i_m}$, and $|b|=b_1+\cdots b_m$. 

Using (3.9) and the fact that $v({\bf i})=F_{i_m}^{a_m}
\cdot v({\bf i}_{m-1})$, we rewrite (3.8) as follows:
$$cv=E_{i_1}^{b_1}\cdots E_{i_m}^{b_m}F_{i_m}^{a_m}\cdot 
v({\bf i}_{m-1})=\sum F_{i_m}^{|b'|-|b|+a_m}E_{i_1}^{b'_1}
\cdots E_{i_m}^{b'_m}p_{b'}
\cdot v({\bf i}_{m-1})\eqno (3.10)$$ 
where the sum is over all $(b'_1,\ldots,b'_m)$ such that $b'_k=b_k$ 
whenever $i_k\ne i_m$. 

It follows that, for some $b'$, we have  
$$F_i^{|b'|-|b|+a_m}E_{i_1}^{b'_1}\cdots E_{i_{m-1}}^{b'_m}p_{b'}
\cdot v({\bf i}_{m-1})=c'v$$
with $c'\in k(q), c'\ne 0$.  
By (3.7), we have $|b'|-|b|+a_m=0$ and $b'_m=0$; also,  
$$E_{i_1}^{b'_1}\cdots E_{i_{m-1}}^{b'_{m-1}}\cdot 
v({\bf i}_{m-1})=c''v$$ 
with $c''\ne 0$. 
Remembering  the inductive assumption, we see that $b'_k=a_k$ for 
all $k\le m-1$. Thus, $b_k=b_k'=a_k$ for all $k\le m-1$ such that 
$i_k\ne i_m$. This completes {\bf Step 1}.

\smallskip

\noindent {\bf Step 2}. Let us prove that $b_m=a_m$. 
If $i_k \ne i_m$ for $k =1, \ldots, m-1$  then the equality 
$b_m=a_m$ follows by comparing degrees.   
So we can assume that $i_k=i_m$ for some $k<m$.
Let $k<m$ be the maximal index such that $i_k=i_m$. Clearly, $k\le m-2$
since ${\bf i}$ is reduced.  By {\bf Step 1}, we have 
$b_{k+1}=a_{k+1},~b_{k+2}=a_{k+2},\ldots, 
b_{m-1}=a_{m-1}$. Combining this observation with (3.7), we can rewrite 
the left hand 
side of
(3.10) as follows. 
$$cv=dE_{i_1}^{b_1}\cdots E_{i_{m-1}}^{b_{m-1}}F_{i_m}^{a_m-b_m}\cdot 
v({\bf i}_{m-1})=dE_{i_1}^{b_1}\cdots E_{i_k}^{b_k} 
E_{i_{k+1}}^{a_{k+1}}\cdots E_{i_{m-1}}^{a_{m-1}}F_{i_m}^{a_m-b_m}
\cdot v({\bf i}_{m-1})$$
for some $d\in k(q),d\ne 0$. 

Then, by  the commutativity property $E_{i_l}F_{i_m}=F_{i_m}E_{i_l}$ for
$k<l<m$, the previous expression is equal to
$$cv=dE_{i_1}^{b_1}\cdots E_{i_k}^{b_k} 
F_{i_m}^{a_m-b_m}E_{i_{k+1}}^{a_{k+1}}\cdots E_{i_{m-1}}^{a_{m-1}}
\cdot v({\bf i}_{m-1})=d'E_{i_1}^{b_1}\cdots E_{i_k}^{b_k}
F_{i_m}^{a_m-b_m}\cdot v({\bf i}_k)\eqno (3.11)$$
(we again used the property (3.7) of the extremal vectors). Since
$i_k=i_m$, it follows that $F_{i_m}\cdot v({\bf i}_k)=0$.  
Hence, the right hand side of  (3.11) is zero unless $a_m-b_m=0$. 
This completes {\bf Step 2}. 

\smallskip

\noindent {\bf Step 3}. Now we are able to complete the proof. 
Since $b_m=a_m$, (3.7) and (3.11) imply that 
$$c_mE_{i_1}^{b_1}\cdots E_{i_{m-1}}^{b_{m-1}}\cdot 
v({\bf i}_{m-1})=cv$$
with some nonzero constants $c,c_m$. 
We conclude that $(b_1,\ldots,b_{m-1})=(a_1,\ldots,a_{m-1})$ by 
the inductive assumption. Combining this with {\bf Step 2}, we 
see that $(b_1,\ldots,b_m)=(a_1,\ldots,a_m)$. 

Proposition 3.3 and Theorem 3.1 are proved.   $\triangleleft$

\bigskip
 
\noindent {\bf Appendix. Skew-fields of fractions and skew polynomials}

\medskip 
 
Let $A$ be an associative ring with unit without zero-divisors.
As in [7], A.2, we say that $A$ satisfies the {\it right Ore}
condition if $aA \cap bA \ne \{0\}$ for any non-zero $a, b \in A$.
The set of right fractions ${\cal F}(A)$ is defined as
the set of all pairs  
$(a,b)$ with $a,b\in A, b\ne 0$ 
modulo the following equivalence relation: 
$(a,b)\sim (c,d)$ if there are 
$f,g\in A\setminus \{0\}$ such that $af=cg$ and $bf=dg$. 
The equivalence class of 
$(a,b)$ in ${\cal F}(A)$ is denoted by $ab^{-1}$.
The ring $A$ is naturally embedded into ${\cal F}(A)$ via
$a \mapsto (a,1)$. 
It is well known that if $A$ satisfies the right Ore
condition then the addition and multiplication in $A$ extend to 
${\cal F}(A)$ so that ${\cal F}(A)$ becomes a skew-field.

\medskip 
 
Now we suppose that $A$ is an algebra over a field ${\bf k}$ with an 
increasing filtration $\big({\bf k}= A_0\i A_1 \i \cdots\big)$, where
each $A_k$ is a finite dimensional ${\bf k}$-vector space, 
$A_kA_l\i A_{k+l}$, and $A = \cup A_k$. We say that $A$ has {\it
polynomial growth} if for all $n\ge 0$ we have $\dim A_n \le p(n)$,
where $p(x)$ is a polynomial.
For the convenience of the reader, we will present a proof of the
following well known lemma (see, e.g., [5]).  

\proclaim Lemma A1. 
Any algebra of polynomial growth without zero-divisors satisfies the right
Ore condition.

\noindent {\sl Proof.}
Assume, on the contrary, that $aA \cap bA = \{0\}$ for some non-zero 
$a, b \in A$.
Denote $I_n=I\cap A_n$ for any subspace $I\i A$. Choose 
some $k$ such that $a,b\in A_k$. 
Then $(aA)_{n+k}\supset a A_n$ and  $(bA)_{n+k}\supset b A_n$, 
which implies
$$\dim (aA)_{n+k}\ge \dim A_n, \,\, \dim (bA)_{n+k}\ge \dim A_n \ .$$
On the other hand, since $aA \cap bA = \{0\}$, it follows that
$$\dim A_{n+k}\ge \dim (aA)_{n+k} + \dim (bA)_{n+k}\ge 2\dim A_n$$
for all $n$.  
Iterating this inequality, we see that $\dim A_{mk}\ge 2^m$ for $m \ge 0$.
This contradicts the condition that $A$ has polynomial growth.
Lemma A1 is proved.   $\triangleleft$

\medskip

Lemma A1 implies that any subalgebra of an algebra $A$ of polynomial 
growth without zero-divisors also satisfies the right
Ore condition.

In particular, consider the ${\bf k}$-algebra $P$ of skew 
polynomials generated by $t_1,\ldots,t_m$ subject to the relations
$t_lt_k=q_{kl}t_lt_k$ for $1\le k<l\le m$, 
where the $q_{kl}$ are some non-zero elements of ${\bf k}$. 
It is easy to see that $P$ has no zero-divisors and has 
polynomial growth with respect to the filtration 
$({\bf k}=P_0\i P_1 \i \cdots )$, where $P_n$ 
is the linear span of 
all monomials in $t_1, \ldots, t_m$ of degree $\le n$.
We see that every subalgebra ${\cal B}$ of $P$  satisfies the right
Ore condition.
Therefore, ${\cal F}({\cal B})$ is a skew subfield of ${\cal F}(P)$.

\bigskip

\centerline {\bf References}

\medskip

\item {[1]}  A.Berenstein, S.Fomin, A.Zelevinsky, Parametrizations 
of canonical bases and totally positive matrices. 
 To appear in {\sl Adv. in Math.} 

\item {[2]} A.Berenstein, A.Zelevinsky, String bases for quantum groups
of type $A_r$,  {\sl Advances in Soviet Math.}, {\bf 16}, 
Part 1 (1993),
51--89.

\item {[3]} A.Berenstein, A.Zelevinsky, Canonical bases for the 
quantum group of type $A_r$ and piecewise-linear combinatorics. 
To appear in {\sl Duke Math. J.}.

\item {[4]} L.D.Faddeev, N.Yu.Reshetikhin, L A.Takhtadzhyan,  
Quantization of Lie groups and
Lie algebras. (Russian) {\sl Algebra i Analiz}, {\bf 1} (1989), 
no. 1, 178--206.  

\item {[5]} K.Iohara, F.Malikov,
Rings of skew polynomials and Gel'fand-Kirillov conjecture for quantum 
groups, {\sl Commun. Math. Phys.} {\bf 164} (1994), 217--238.

\item {[6]} A.Joseph, Sur une conjecture de Feigin.
 {\sl C. R. Acad. Sci. Paris S\'er.I Math.} {\bf 320} (1995), no. 12, 1441--1444.

\item {[7]} A.Joseph, Quantum groups and their primitive ideals, 
Springer-Verlag, Berlin, 1995. 

\item {[8]} G.Lusztig, 
Introduction to quantum groups,
Birkhauser, Boston, 1993. 

\item {[9]} G.Lusztig, Problems on canonical bases. Algebraic groups and 
their generalizations: quantum and infinite-dimensional methods 
(University Park, PA, 1991),
169--176, {\sl Proc. Sympos. Pure Math.}, {\bf 56}, Part 2, 
Amer. Math. Soc., 
Providence, RI, 1994.

\item {[10]} S.Majid, Algebras and Hopf algebras in braided 
categories. Advances in Hopf
algebras (Chicago, IL, 1992), 55--105, {\sl Lecture Notes in 
Pure and Appl. Math.}, {\bf 158}, Dekker, New York,
1994.

\item {[11]} A.Morozov, L.Vinet,  Free-field representation of group 
element for simple quantum
group. Preprint ITEP-M3/94, CRM-2202, {\it hep-th} 9409093.   

\item {[12]}  A.Panov, Skew fields of twisted rational functions and 
the skew field of rational functions on ${\rm GL}\sb q(n,K)$. {\sl St.
Petersburg Math. J}, {\bf 7} (1996), no. 1.

\end